\documentclass[a4paper, 10pt ]{article}

\usepackage[a4paper,left=3cm,right=3cm,top=2.5cm,bottom=2.5cm]{geometry}
\usepackage{hyperref}
\hypersetup{colorlinks,linkcolor={blue},citecolor={blue},urlcolor={black}}
\usepackage{subfig}
\usepackage{pgfplots}
\usepackage{pgfplotstable}
\usepackage{mathtools}
\usepackage{multicol}
\usepackage{comment}
\usepackage{booktabs}
\pgfplotsset{compat=1.5}
\usepackage{amssymb}
\usepackage{url}
\usepackage{bm}

\usepackage{pdfsync}
\usepackage{float}
\usepackage{tabularx}
\usepackage{enumerate}
\usepackage{array}
\usepackage{xspace}
\usepackage{tikz}
\usepackage{tikz-cd}
\usepackage{tikzsymbols}
\usetikzlibrary{calc,trees,positioning,arrows,chains,shapes.geometric,%
    decorations.pathreplacing,decorations.pathmorphing,shapes,%
    matrix,shapes.symbols, decorations.markings, patterns,fit}

\usepackage{siunitx}

\usepackage{authblk}

\usepackage[draft,inline,marginclue]{fixme}
\usepackage{mathrsfs}
\usepackage{color, colortbl}
\usepackage{multirow}

\usepackage{amsthm}
\usepackage{amsmath}
\usepackage{stmaryrd}
\usepackage[ruled,vlined,linesnumbered]{algorithm2e}
\usepackage{graphicx}
\usepackage{booktabs}
\usepackage{cleveref}

\FXRegisterAuthor{gs}{ags}{\color{red}GS}
\FXRegisterAuthor{fr}{afr}{\color{blue}FR}

\newtheorem{theorem}{Theorem}

\theoremstyle{remark}
\newtheorem{rmk}{Remark}

\newcommand{\x}{\ensuremath{\mathbf{x}}}

\DeclareMathOperator*{\argmin}{\arg\!\min}
\DeclareMathOperator*{\argmax}{\arg\!\max}

\newcolumntype{C}[1]{>{\centering\arraybackslash}m{#1}}

\definecolor{Gray}{gray}{0.9}

\begin{document}

\title{Non-linear manifold ROM with Convolutional Autoencoders and Reduced Over-Collocation method}

\author[]{Francesco~Romor\footnote{francesco.romor@sissa.it}}
\author[]{Giovanni~Stabile\footnote{giovanni.stabile@sissa.it}}
\author[]{Gianluigi~Rozza\footnote{gianluigi.rozza@sissa.it}}

\affil{Mathematics Area, mathLab, SISSA, via Bonomea 265, I-34136 Trieste,
  Italy}

\maketitle

\begin{abstract}
  Non-affine parametric dependencies, nonlinearities and advection-dominated
  regimes of the model of interest can result in a slow Kolmogorov n-width
  decay, which precludes the realization of efficient reduced-order models based
  on linear subspace approximations. Among the possible solutions, there are
  purely data-driven methods that leverage autoencoders and their variants to
  learn a latent representation of the dynamical system, and then evolve it in
  time with another architecture. Despite their success in many applications
  where standard linear techniques fail, more has to be done to increase the
  interpretability of the results, especially outside the training range and not
  in regimes characterized by an abundance of data. Not to mention that none of
  the knowledge on the physics of the model is exploited during the predictive
  phase. In order to overcome these weaknesses, we implement the non-linear
  manifold method introduced by Carlberg et al~\cite{lee2020model} with hyper-reduction achieved
  through reduced over-collocation and teacher-student training of a reduced
  decoder. We test the methodology on a 2d non-linear conservation law and
  a 2d shallow water models, and compare the
  results obtained with a purely data-driven method for which the
  dynamics is evolved in time with a long-short term memory network.
\end{abstract}

\tableofcontents
\listoffixmes

\section{Introduction}
\label{sec:intro}

In real world engineering scenarios, when performing
outer loop applications such as optimization, uncertainty quantification,
sensitivity analysis, parametric partial differential equations (PDEs) often need to be solved numerically numerous times. However, relying on the mathematical properties of some
parametric PDEs, the computational cost for many query problems can be drastically
reduced taking into account previous results on a set of training parameters:
the procedure for the design of reduced-order models (MORs) is divided in an
offline (training) stage, during which a set of training solutions is collected,
and an online (testing or predictive) stage, which employs the compressed
information from the previous step to predict the solutions of the PDE of
interest for unseen parameters.
This reduction is performed numerically defining a low-dimensional global basis devised
in the offline stage, and can be carried out independently from the class of
numerical methods chosen: finite element (FEM), spectral element (SEM),
discontinuous Galerkin (DGM), and finite volumes method (FVM). One of the most
employed model-order reduction method (MOR) is the reduced basis
method~\cite{hesthaven2016certified, quarteroni2015reduced}.

Depending on the parametric dependency and mathematical nature of some PDEs,
various issues may occur: the Kolmogorov n-width (KnW) is used to characterize the
approximability of the solution manifold, that is the set of parameter-dependent
solutions of the PDE, by a linear trial subspace. A slow decaying KnW is a
symptom of the difficulties in the design of efficient ROMs: this results in
the necessity of using a high number of reduced basis or proper orthogonal
decomposition (POD) method's modes,
corrupting the efficiency of the ROMs till the point that the gain into the computational cost becomes irrelevant. One class of PDEs where this
behaviour is evident are time-dependent advection-dominated PDEs. Moreover, non-linear
PDEs require hyper-reduction procedures to make the reduced equations
independent from the number of degrees of freedom of the full-order model (FOM).

Recently, leveraging machine learning's advances in manifold learning, a class
of ROMs that employ a non-linear trial manifold built with
convolutional autoencoders (CAEs)~\cite{goodfellow2016deep} was developed by Carlberg et
al~\cite{lee2020model}. One of the benefits of this approach is the possibility
to employ a small latent dimension of the ROMs,
thus overcoming the slow decay of the KnW for some parametric PDEs, at the
expense of introducing additional non-linearities from the neural networks (NNs)
and sometimes more substantial training costs in the offline stage. The properties of the
non-linear manifold methods include the need of less stabilization mechanisms,
the less intrusiveness on the FOM solvers --- they are in fact equations-based
rather than fully-intrusive --- and the possibility to apply them for a much
broader class of parametric PDEs, differently from ROMs devised specifically for advection-dominated
problems.

A hyper-reduction scheme for non-linear manifold Least-Squares Petrov-Galerkin
(NM-LSPG) and non-linear manifold Galerkin (NM-G), is introduced
in~\cite{kim2020fast}: it relies on Gauss-Newton tensor approximation (GNAT)~\cite{Carlberg_2013}
hyper-reduction method and shallow masked autoencoders to select only
the degrees of freedom that explain the dynamics and therefore restrict in an
efficient way the decoder and the discretized residuals. As we will see in our test cases, the
reconstruction error of the autoencoder employed empirically bounds from below
the errors of all the other non-linear manifold ROMs built upon. Therefore, we devise a
method that is independent on the choice of the architecture: a sparse shallow
autoencoder is not necessary anymore, and any NN architecture, like CAEs, could be in
principle employed. This frees the way to imposing additional inductive biases
that help to speed-up the offline stage and to achieve
accurate approximations of the discrete solution manifolds, a crucial
requirement. Moreover, in some cases, reconstructing the residuals with GNAT is
not efficient, still because of a slow decaying KnW, so we choose to employ the reduced
over-collocation hyper-reduction method (ROC)~\cite{chen2021eim}: in this case
the equation's numerical residual is not reconstructed on a global basis, but
collocated on some nodes of the mesh called collocation points.

Once the CAE reaches a satisfactory approximation of the discrete solution
manifold, purely data-driven NN PDE can be trained to predict the
latent dynamics: the gold standard that is being established for this task are
long-short term memory networks (LSTM). Their online computational cost is low
even w.r.t linear ROMs, but some new issues appear: their accuracy depends on
the regularity of the latent dynamics, especially when predicting the solutions
for parameters outside the training range, in the extrapolation regimes; they
require hyper-parameters tuning, and all the connections to the PDEs model are
completely lost, resulting in a loss of interpretability of the results. The
non-linear manifold hyper-reduced ROMs we develop solve these issues, at the
expense of a higher computational cost in the online stage, since at each time
step a physics-based residual is minimized. Moreover, a
posteriori error estimates are available~\cite{lee2020model, kim2020fast}. In
our test cases, we compare these two approaches to enlight their differences,
weak and strong points.\\

The structure of this paper is as follows. In section~\ref{sec:ml} we delve
into the topic of manifold learning which has many connections with
reduced-order modelling, especially since the recent entry of machine
learning in the design of ROMs. We will proceed introducing the Kolmogorov
n-with (KnW), and we will show that some classes of parametric PDEs suffer from the so
called slow decaying Kolmogorov n-width.
In
section~\ref{sec:nmlspgroc}, the non-linear manifold (NM) reduced-order models
based on the work of Carlberg et al~\cite{lee2020model} are introduced. We will
focus on the non-linear manifold least-squares Petrov-Galerkin method (NM-LSPG).
Afterwords, we describe our new hyper-reduced ROMs: NM-LSPG with reduced
over-collocation (NM-LSPG-ROC) and NM-LSPG with reduced over-collocation and
teacher-student training of a compressed decoder (NM-LSPG-ROC-TS).
These two approaches, to the best of the authors' knowledge, are introduced here 
for the first time.
 In
section~\ref{sec:results}, the new model order reduction (MOR) methods are tested on a 2d parametric
non-linear time-dependent conservation law model and a 2d parametric non-linear
time-dependent shallow water equations model. In section~\ref{sec:discussion}
follows a discussion on the results obtained and in the conclusive
section~\ref{sec:conclusions} some possible future directions of research are explored.

\section{Manifold learning}
\label{sec:ml}
The subject of manifold learning, classified as a topic of machine learning,
had its unique flavour in model order reduction even before nowadays breakout of
scientific machine learning~\cite{baker2019workshop}. The workhorse of the model order reduction
community is POD or SVD. A lot of real applications though,
required new methods to approximate the solution manifold in a non-linear
fashion. The symptom of this behaviour is a slow decaying Kolmogorov n-width. Some approaches rely on the locality of the validity region of a linear
approximation with POD, both in the parameter space and in the spatial and
temporal domains, others implement non-linear dimensionality reduction methods from
machine learning~\cite{murphy2012machine}: kernel principal component analysis (KPCA), Isomap, clustering algorithms. Non-linear MORs include approximations by
rational functions, splines or other non-linear functions collected in a dictionary~\cite{benner2017model}.

Interpolatory approaches of the solution manifold with respect to the parameters
have been developed, sometimes combined with non-linear dimension reduction
techniques like KPCA and its variants:
interpolation with geodesics on the Grassmann
manifold~\cite{amsallem2008interpolation}, interpolation on the latent space
obtained with Isomap dimensionality reduction
method~\cite{franz2014interpolation, bhattacharjee2016nonlinear}, interpolation
with optimal transport~\cite{bernard2018reduced}, dictionary-based ROM that make
use of clustering in the Grassmmannian manifold and classification with neural networks
(NN)~\cite{diez2021nonlinear}, local kernel principal component
analysis~\cite{li2021model}. At the same time, domain decomposition approaches
tackled locality in space~\cite{lucia2003reduced, buffoni2009iterative}.

One particular class of dimension reduction techniques is represented by
autoencoders, and more generally by other architectures that rely on NNs. In the
recent literature many achievements are brought by CAEs, and by extension Generative Adversarial Networks (GANS),
Variational Autoencoders (VAEs), bayesian convolutional autoencoders~\cite{goodfellow2016deep}: in~\cite{mucke2021reduced} convolutional autoencoders are utilized for
dimensionality reduction and long-short Term Memory (LSTM) NNs or causal
convolutional neural networks
are used for time-stepping; in~\cite{fresca2022pod} the evolution of the
dynamics and the parameter dependency is learned at the same time of the latent
space with a forward NN and a CNNs on randomized SVD compressed snapshots,
respectively; and in~\cite{xu2020multi} spatial and temporal features are
separately learned with a multi-level convolutional autoencoder.

In order to extend this architectures to datasets not structured in orthogonal
grids, geometric deep learning~\cite{bronstein2017geometric} is called to the task. There are not many works on geometric deep learning applied to model
order
reduction for mesh-based simulations that achieve the same accuracy of CNNs,
yet. Promising results are reached by an architecture that employs graph neural
networks (GNNs) and a physics-informed loss~\cite{pfaff2020learning}. In the
future, the potential of GNNs will probably be leveraged extending the range of
applicability of nowadays frameworks.

The setting we will base our studies on, does not depend directly on the
numerical method used to discretize the parametric PDEs at hand (FVM, FEM, SEM, DGM), so the mathematical formulation will generically be founded on models
represented by a parametric system of time-dependent (but also time-independent)
PDEs, consisting of a non-linear parametric differential operator $\mathcal{G}$
and of the boundary differential operators $\mathcal{B}, \mathcal{B}_{0}$ that
represent the boundary and initial conditions respectively,
\begin{align}
    \forall\boldsymbol{\mu}\in\mathcal{P}\qquad
    \begin{cases}
        \mathcal{G}(\boldsymbol{\mu}, U(\boldsymbol{\mu}, t, \mathbf{x}))&=0 \qquad (\mathbf{x}, t)\in\Omega\times [0, T], \\
        \mathcal{B}(\boldsymbol{\mu}, U(\boldsymbol{\mu}, t, \mathbf{y}))&=0 \qquad (\mathbf{y}, t)\in\partial\Omega\times[0, T],\\
        \mathcal{B}_{0}(\boldsymbol{\mu}, U(\boldsymbol{\mu}, 0, \mathbf{x}))&=0 \qquad (\mathbf{x}, t)\in\Omega,
    \end{cases}
    \label{eq:PDEsystem}
\end{align}
where $\mathcal{P}$ is the parameter space, $U$ are the state variables and
$\Omega$ is the 2D
or 3D spatial domain. This formulation includes also
coupled systems of PDEs. We assume that the solutions belong to a certain Banach
or Hilbert space $Y$, varying $(\boldsymbol{\mu}, t)\in \mathcal{P}\times [0,
T]$. The solution manifold $\mathcal{M}$ is represented by the set
\begin{equation}
    \mathcal{M} = \{U(\boldsymbol{\mu}, t)\in Y \ | \ \boldsymbol{\mu}\in\mathcal{P}, t\in [0, T]\}.
\end{equation}
\subsection{Approximability by n-dimensional subspaces and Kolmogorov n-width}

We want to remark some results available in the literature, in order to state and
comment, for our needs, the problem of solution manifold approximability. In
particular, we will prove that for a certain class of parametric PDEs, whom our
benchmarks belong to, the Kolmogorov n-width decays slowly. Thus, classical
Petrov-Galerkin projection with POD needs to be overcome with non-linear
methods in place of POD to achieve efficient
ROMs.

Let $(X, \lVert\cdot\rVert)_{X}$ be a complex Banach space, and $K\subset X$ a
    compact subspace, the Kolmogorov n-width (KnW) of $K$ in $X$ is defined as
\begin{equation}
    d_{n}(K)_{X}=\inf_{\text{dim}(W)=n}\max_{v\in K}\min_{w\in W}\lVert v-w \rVert_{X}.
\end{equation}
Let $(Y, \lVert\cdot\rVert_{Y})$ be a complex Banach space and
$L:K\subset\subset X \rightarrow Y$. In our framework $K$ is the parameter
space, possibly infinite dimensional and $L$ is the solution map of the system
of parametric PDEs at hand, from the parameter space to the solution manifold.
In order to define $L$ we have to suppose that for each parameter in
$K$ there is a unique solution in $Y$.

Following \cite{cohen2016kolmogorov}, it can be proved that for holomorphic $L$,
thus not necessarily linear, the Kolmogorov n-width decay is one polynomial
order below the Kolmogorov n-width of the parameter space $K$.

\begin{theorem}[Theorem 1, \cite{cohen2016kolmogorov}]
    \label{theo:analytic solution map}
    Suppose u is a holomorphic mapping from an open set $O\subset X$ into $Y$
    and u is uniformly bounded on $O$,
    \begin{equation}
        \exists B>0,\quad\sup_{\x\in O} \lVert u(x)\rVert_{Y} \leq B.
    \end{equation}
    If $K\subset O$ is a compact subset of $X$, then for any $s>1$ and $t<s-1$,
    \begin{equation}
        \sup_{n\geq 1} n^{s} d_{n}(K)_{X} < \infty \Rightarrow \sup_{n\geq 1} n^{t}d_{n}(u(K))_{K} < \infty.
    \end{equation}
\end{theorem}

In particular, if the hypothesis of the previous theorem are satisfied and if
$K$ is a finite dimensional linear subspace, the Kolmogorov n-width decay is
exponential. In general, elliptic PDEs, affinely decomposable with respect to the
parameters, satisfy the hypothesis of the previous
theorem~\cite{babuvska2007stochastic, lassila2013generalized}. Not always
nonlinearities cause a slow decaying KnW: using theorem~\ref{theo:analytic
solution map}, in~\cite{cohen2016kolmogorov} they prove that the parametric PDE
on a bounded Lipschitz domain $\Omega\subset\mathbb{R}^{3}$
\begin{equation}
    u^{3}-\nabla \cdot(\exp{a})\nabla u = f,\quad K=\{a\in L^{\infty}(\Omega)|\lVert a\rVert_{C^{\alpha}}\leq M\},\quad f\in H^{-1}(\Omega),
\end{equation}
with homogeneous Dirichlet boundary conditions, where $C^{\alpha}$ is the space of
H\"older functions, satisfies the hypothesis of the previous theorem. Actually,
for H\"older functions the KnW is bounded above by $n^{-\alpha/3}$,
which is not a fast convergence,
but if instead $a$ belongs to the Sobolev space $W^{s, \infty}(\Omega)$ then the
upper bound is $n^{s/3}$~\cite{geller2014kolmogorov}. We will consider a good KnW
decay if it has a higher infinitesimal order than $n^{-1}$.

The same is not valid in the case of the simplest linear advection problem. We
briefly report some results from the literature on classical hyperbolic PDEs,
for $(t,\mu)\in K=[0,1]^2$ with the standard norm and $Y=L^{2}([0, 1])$,
\begin{align}
    L:(t,\mu) \mapsto u(t, x, \mu) \quad\text{s.t.}\quad\partial_{t}u - \mu \partial_{x}u = 0\qquad           & d_{n}(\mathcal{M})_{L^{2}}> n^{-\frac{1}{2}},\ \text{from}~\cite{ohlberger2015reduced},        \\
    L:(t,\mu) \mapsto u(t, x, \mu) \quad\text{s.t.}\quad\partial^{2}_{tt}u - \mu \partial^{2}_{xx}u = 0\qquad & d_{n}(\mathcal{M})_{L^{2}}> n^{-\frac{1}{2}}, \ \text{from}~\cite{greif2019decay}.
\end{align}

This behaviour is not restricted only to advection-dominated problems.
Intuitively also solution manifolds that are characterized by a parametrized locality in space suffer from slow decaying KnW, like elliptic problems
with singular sources parametrically moving in the domain~\cite{franco2021deep}.

Our newly developed ROMs, should solve the issue of slow decaying KnW in the applications,
guaranteeing a low latent or reduced dimension of the approximate solution
manifold. This, because the linear trial manifold, frequently generated by the
leading POD modes, is
substituted with a non-linear trial manifold, parametrized by the decoder of an
autoencoder. The test cases we present in section~\ref{sec:results}, were chosen
in order to be advection-dominated and particularly not suited for linear ROMs,
as shown in~\cite{kim2020fast} for the 2d Burgers' equation, that has a
close relationship with the NCL problem.

\begin{rmk}[Extensions of the Kolmogorov n-width to non-linear approximations]
    The autoencoders we implement in our test cases are at least continuous as
    composition of continuous activations and linear functions. In the
    literature there exists possible non-linear extensions of the KnW such as the manifold
    width~\cite{devore1989optimal},
    \begin{equation}
        \delta_{n}(L(K), X) = \inf_{\substack{\psi\in\mathcal{C}(X, \mathbb{R}^{n})\\\phi\in\mathcal{C}(\mathbb{R}^{n}, X)}}\sup_{x\in L(K)}\lVert x-(\phi\circ\psi)(x)\rVert_{X},
    \end{equation}
    library widths~\cite{temlyakov1998nonlinear}, and entropy numbers, for more
    insights see~\cite{benner2017model}.
\end{rmk}

\subsection{Singular values decomposition and discrete spectral decay}
From the discrete point of view the same problematic in tackling the
reduction of parametric PDEs with slow KnW decay is encountered: in this
case the discrete solution manifold is actually the set of discrete solutions of the
full-order model for a selected finite set of parameter values.
Singular Value Decomposition (SVD) or eigenvalue decompositions (for
symmetric positive definite matrices), usually employed on the snapshots
matrix for the evaluation of the reduced modes, are characterized by the
fact that modes are linear combinations of the snapshots. This is not enough
to approximate snapshots that are orthogonal with respect to the collection
of the training set.

In practice, looking at the residual energy retained by the discarded modes,
is an indicator of approximability with linear subspaces. Let us assume that
$d$ is the total number of degrees of freedom
of the discrete problem, $N_{\text{train}}$
is the number of training snapshots, and, $r$, the reduced
dimension such that, $r\leq N_{\text{train}}\ll d$. Then,
$X\in\mathbb{R}^{d}\times\mathbb{R}^{N_{\text{train}}}$ is the matrix that has the
training snapshots $\{U_i\}_{i=1}^{N_{\text{train}}}$ as columns,
$V\in\mathbb{R}^{d}\times\mathbb{R}^{r}$ are the reduced modes from the SVD
of $X$, and
$\{\sigma_{i}\}_{i=1}^{d}$ are the increasingly ordered singular values. It is valid the following relationship of the residul
energy (to the left) with the KnW (to the right),
\begin{equation}
    \sqrt{\sum_{i=r+1}^{d}\sigma_{i}^{2}}=\lVert (I-VV^{T})X\rVert_{F}\geq \max_{i=1,\dots,N_{\text{train}}}\lVert (I-VV^{T})U_i\rVert_{2}\geq d_{n}(\{U_i\}_{i=1}^{N_{\text{train}}})_{\mathbb{R}^{d}},
\end{equation}
where $\lVert\cdot\rVert$ is the Frobenious norm.

Even though some problems have a slow decaying KnW, this affects only the
asymptotic convergence of the ROMs w.r.t. the reduced dimension, while for
some applications a satisfying accuracy of the projection error is reached
with less than $100$ modes, as shown in our benchmarks. So what is actually
lost in MOR for problems with a slow decaying KnW is the fast convergence of
the projection error w.r.t. the reduced dimension, not the possibility to
perform a MOR with enough modes. Moreover, the discrete solution manifold's
KnW depends on the time step and spatial discretization size, so that, especially
when a coarse mesh is employed, the Knw decays faster w.r.t. the KnW of the
continuous solution manifold.
\subsection{Convolutional Autoencoders}
We have chosen to overcome the slowly decaying KnW problem employing
autoencoders~\cite{goodfellow2016deep} as non-linear dimension reduction method
substituting POD. Some ROMs predict the latent dynamics on a linear trial
manifold with artificial neural networks~\cite{PiBaRoHe2021}, so they are still
classified as linear ROMs and, in fact, they are still affected by the slow KnW
decay. Moreover, we remark that while some non-linear approaches to model order
reduction are specifically tailored for advection-dominated
problems~\cite{torlo2020model, PapapiccoDemoGirfoglioStabileRozza2021},
autoencoders are a more general approach. On the other hand, they are also
particularly suited to advection-dominated problems with respect to local and/or
partitioned ROMs that implement domain decomposition, even when non-linear dimension reduction techniques are employed
locally~\cite{li2021model}: this is because considering a non discrete
parametric space, like the time interval of a simple linear advection problem, an infinite number of local linear and/or non-linear ROMs would be needed to counter
the slow decaying KnW. As anticipated, we implement convolutional
autoencoders~\cite{goodfellow2016deep} in libtorch the C++ frontend of
PyTorch~\cite{NEURIPS2019_9015}.

Let us define $X_h\subset\mathbb{R}^{d}$ as the state discretization space with
$d$ the number of degrees of freedom and $h$ the discretization step. The
snapshots are divided in a training set
$\mathcal{U}_{\text{train}}=\{\mathbf{U}_{i}\}_{i=1,\dots,N_{\text{train}}}\subset
X_h$ and a test set
$\mathcal{U}_{\text{test}}=\{\mathbf{U}_{i}\}_{i=1,\dots,N_{\text{test}}}\subset
X_h$. If the problem has different states and/or vectorial states the training
and test set are split in channels for each state and/or component. For example
in the 2d non-linear conservation law, we consider two channels, one for each velocity
component. In the shallow water test case, we consider three channels, one for
each velocity component and one for the free surface height. So, in general, we reshape the snapshots such that
$\mathcal{U}_{\text{train}},\ \mathcal{U}_{\text{test}}\subset
(\mathbb{R}^{d/c})^c$ where $c$ is the number of channels.

As preprocessing step the snapshots are centered and normalized to assume values
in the interval $[-1, 1]$
\begin{equation}
    \frac{2}{\mathcal{U}_{\text{max}}-\mathcal{U}_{\text{min}}}\left(\mathcal{U}_{\text{train}}-\mathcal{U}_{\text{mean}} - \frac{\mathcal{U}_{\text{min}}+\mathcal{U}_{\text{max}}}{2}\right),
\end{equation}
where $\mathcal{U}_{\text{mean}},\ \mathcal{U}_{\text{max}},\
\mathcal{U}_{\text{min}}\in(\mathbb{R}^{d/c})^c$ are evaluated channel wise. The
same values obtained from the training snapshots $\mathcal{U}_{\text{train}}$, are employed to center and normalize the test set
$\mathcal{U}_{\text{test}}$.

\begin{figure}[h]
    \centering
    \includegraphics[width=0.6\textwidth]{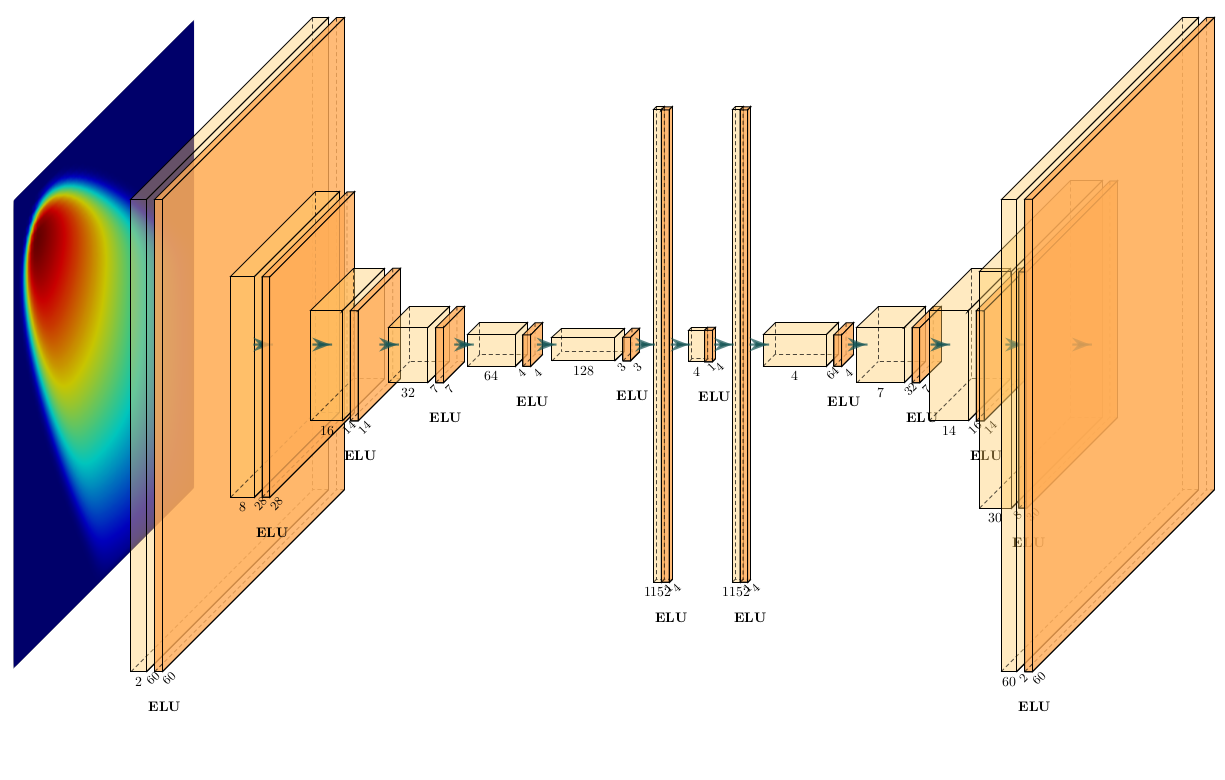}
    \caption{Prototype a convolutional autoencoder architecture employed. The number of layers is the same of the CAE used for the 2d non-linear conservation law test case, but the filter sizes are rescaled for better rendering.}
    \label{fig:cae}
\end{figure}

We define the encoder $\psi:(\mathbb{R}^{d/c})^c\rightarrow\mathbb{R}^{r}$ and
the decoder $\phi:\mathbb{R}^{r}\rightarrow(\mathbb{R}^{d/c})^{c}$, where $r\ll
d$ is the reduced or latent dimension, as neural networks made by subsequent
convolutional layers and linear layers at the end and at the beginning,
respectively. For the particular architecture used in the applications we defer
it to the appendix~\ref{sec:NN}. In Figure~\ref{fig:cae} is represented the
convolutional autoencoder applied for the 2d non-linear conservation law test case, with an
approximate size of the filters and the actual number of layers~\footnote{Figure~\ref{fig:cae} and Figure~\ref{fig:teacher student net}
were made with the open-source package from GitHub~\url{https://github.com/HarisIqbal88/PlotNeuralNet}.}.
\begin{rmk}[Regularity]
    Regarding the regularity of the CAE, related to the choice of activation
    functions, it is proved in Theorem 4.2 from ~\cite{lee2020model} that
    NM-LSTM and non-linear manifold Galerkin methods are asymptotically
    equivalent provided that the decoder is twice differentiable. Since we are
    only employing the NM-LSTM method, our only
    concern in the choice of activation functions is that the reconstruction
    error is sufficiently low, so that the accuracy of the whole procedure is
    not undermined.
\end{rmk}
For each batch $\{\mathbf{U}_{i}\}_{i=1}^b\subset\mathcal{U}_{\text{train}}$,
the loss employed is the sum of the reconstruction error, and a regularizing term
for the weights:
\begin{equation}
    \label{eq:loss cae}
    \mathcal{L}(\{\mathbf{U}_{i}^{t_i}\}_{i=1}^b;\Theta) = \frac{1}{b}\sum_{i=1}^{b}\frac{\lVert \mathbf{U}^{t_i}_{i} - (\phi\circ \psi)(\mathbf{U}^{t_i}_{i})\rVert^{2}_{2}}{\lVert\mathbf{U}^{t_i}_i\rVert^{2}_{2}} + \lambda_{1} \lVert \Theta \rVert^{2}_{2},
\end{equation}
where $\Theta$ represents the weights of the convolutional autoencoder.  The choice of the relative
mean squared reconstruction error is important when, varying the parameter
$\{\boldsymbol{\mu}_i\}_{i=1}^{N_{\text{train}}}$, the snapshots
$\{\mathbf{U}_i\}_{i=1}^{N_{\text{train}}}$ have different orders of magnitude:
for example
this is the case of flows propagating from a local source on the whole domain
with a constant zero state as initial condition.

\begin{figure}[htp!]
    \centering
    \includegraphics[width=1\textwidth]{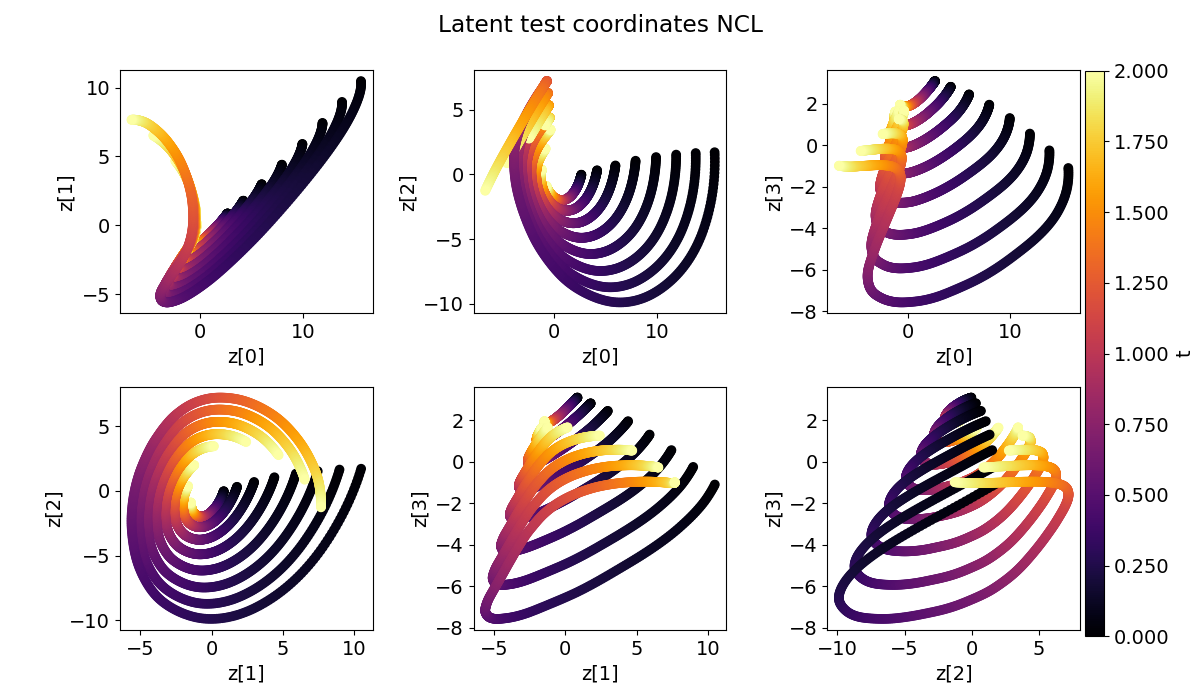}
    \caption{Correlations among the $4$ latent coordinates of the odd parameters of the test set, $8$ in total. They are obtained projecting the test snapshots into the latent space of the CAE with the encoder. The coloring corresponds to the time instants from $0$ to $2$. }
    \label{fig: latent_true_ncl}
\end{figure}

The training is performed with Adam stochastic optimization
method~\cite{kingma2014adam}. After the training the whole evolution of the
dynamics is carried out with a non-linear optimization algorithm minimizing the residual
on the latent domain, as described in section~\ref{sec:nmlspgroc}. For each new
parametric instance and associated initial state $\mathbf{U}_0\in X_h$,
the latent initial condition is found with a single forward of the encoder $z_0
= \psi(\mathbf{U}_0)$ after centering and normalizing $\mathbf{U}_0$. Then, for
each time instant $t$, the decoder
\begin{equation}
    \phi(z(t))=\mathbf{U}(t)\in (\phi\circ\psi)(X_h)\subset\mathbb{R}^{d},
\end{equation}
is used as parametrization of an approximate solution manifold as will be
explained in section~\ref{sec:nmlspgroc}, including in $\phi$ the renormalization of the
output.
\begin{rmk}[Initial condition]
    \label{rmk:initial condition}
    With respect to the initial implementation of Carlberg et
    al~\cite{lee2020model} our approach for the definition of the latent
    parametrization of the solution manifold is different in the management of
    the initial condition. Instead of using directly the decoder $\phi$, they
    define the map from the latent to the state space
    $f:\mathbb{R}^{r}\rightarrow\mathbb{R}^{d}$, including the renormalizations
    in the $\phi$ map
    for brevity, as
    \begin{equation}
        f(z) = \phi(z)+\mathbf{U}_0(\mu)-\phi(\psi(\mathbf{0})),\quad \text{s.t.}\ f(\psi(\mathbf{0}))= \mathbf{U}_0(\mu),
    \end{equation}
    so that the reconstruction is exact at the initial parametric time instances. So, supposing that the initial condition is parametrically dependent, all
    initial conditions coincide to $\psi(\mathbf{0})$ in the latent space, and the decoder has to learn
    variations from the initial latent condition. We prefer instead to
    split the
    initial conditions in the latent space in order to aid the non-linear
    optimization algorithm, and leave to the training of the CAE the accurate
    approximation of the initial condition. This splitting of the initial
    conditions can be seen in Figure~\ref{fig: latent_true_ncl} and \ref{fig: latent_test_swe}. The additional cost of our
    implementation is the forward of the initial parametric condition through
    the encoder as first step.
\end{rmk}
\begin{rmk}[Inductive biases]
    Imposing inductive biases to increase the convergence
    speed of a deep learning model to the desired solution has always been a
    winning strategy in machine learning.
    This is translated in the context of physical models with the  possibility
    to include, among others, the following inductive biases:
    first principles (conservation laws~\cite{lee2019deep}, equations governing the physical phenomenon), geometrical simmetries
    (group invariant filters~\cite{smets2020pde, finzi2020generalizing}),
    numerical schemes/residuals (discrete residuals, latent time advancement
    with Runge Kutta schemes),
    latent regularity (minimize the curvature of latent trajectories), latent
    dynamics (linear or quadratic latent dynamics~\cite{goyal2021lqresnet}). As inductive bias, we will impose the positivity
    of the state variables that are known to be positive throughout their
    trajectory with a final ReLU activation.
\end{rmk}

\section{Evolution of the latent dynamics with NM-LSPG-ROC}
\label{sec:nmlspgroc}
The non-linear manifold method introduced by Carlberg et al~\cite{lee2020model}
does not perform a complete dimension reduction since at each time step the decoder
reconstruct the state from the latent coordinates to the whole domain, still
depending on the number of degrees of freedom of the FOM. We revisit the
non-linear manifold least-squares Petrov-Galerkin method (NM-LSPG) with small modifications
and introduce two novel hyper-reduction procedures: one combines teacher-student
training of a compressed decoder with the reduced over-collocation method
(NM-LSPG-ROC-TS), the other implements only the hyper-reduction of the residual
with reduced over-collocation (NM-LSPG-ROC).
\subsection{Non-linear manifold least-squares Petrov-Galerkin}
\label{subsec: NM-LSPG}
We assume that the numerical method of preference discretizes the
system~\ref{eq:PDEsystem} in space and in time with an implicit scheme, $G_{h,
\delta t}:\mathcal{P}\times X_h\times X_h^{|I_t|}\rightarrow X_h$
\begin{equation}
    G_{h, \delta t}(\boldsymbol{\mu}, \mathbf{U}_h^t, \{\mathbf{U}_h^s\}_{s\in I_t})=\mathbf{0},
    \label{eq:PDEsystemDiscreteU}
\end{equation}
where $h,\ \delta t$ are the spatial and temporal discretization steps chosen,
$X_h\subset\mathbb{R}^{d}$ is the state discretization space ($d$ is the number
of degrees of freedom), and $I_t$ is the set of past state indexes employed in
the temporal numerical scheme to solve for $\mathbf{U}^t_h$. We remark that the
numerical discretization employed can differ from the one used to solve for the
full-order training snapshots. The method is thus equations-based rather than
fully intrusive. This will be the case for the 2d shallow water equations model
in Subsection~\ref{subs:shallow water}.

For each discrete time instant $t$ the following non-linear least-squares problem
is solved for the latent state $\mathbf{z}^t \in Z$, with the Levenberg-Marquardt algorithm~\cite{quarteroni2010matematica}
\begin{equation}
    \mathbf{z}^t = \argmin_{\mathbf{z}\in \mathbb{R}^r}\ \lVert G_{h, \delta t}(\boldsymbol{\mu}, \phi(\mathbf{z}), \{\phi(\mathbf{z}^s)\}_{s\in I_t})\rVert^{2}_{X_h}
    \label{eq:PDEsystemDiscrete}.
\end{equation}
That is for each time instant the following intermediate solutions
$\{\mathbf{z}^{t,k}\}_{k\in\{0,\dots,N(t)\}},\ \mathbf{z}^{t, 0} =
\mathbf{z}^{t-1, N(t-1)}$ of the linear system in $\mathbb{R}^r$ are computed,
\begin{align}
    \left((dG^{t, k-1}d\phi^{t, k-1})^T dG^{t, k-1}d\phi^{t, k-1} + \lambda\ I_d\right) \delta\mathbf{z}^{t,k} = -(dG^{t, k-1}d\phi^{t, k-1})^T G^{t, k}\label{eq:levenbergMarquardt},\\
    \mathbf{z}^{t,k} = \mathbf{z}^{t,k-1} + \alpha^k\ \delta\mathbf{z}^{t,k},
\end{align}
where $d\phi^{t, k-1}=d\phi(\mathbf{z}^{t, k-1})\in\mathbb{R}^{d\times r}$ and
$dG^{t, k-1}=dG_{h, \delta t}(\boldsymbol{\mu}, \phi(\mathbf{z}^{t, k-1}),
\{\phi(\mathbf{z}^{s})\}_{s\in I_t})\in\mathbb{R}^{d\times d}$ are the Jacobian
matrices of the embedding map $\phi$ and of the residual $dG_{h, \delta t}$
respectively, $\lambda$ is a factor that evolves during the non-linear
optimization and balances between a Gauss-Newton and a steepest descent method,
and finally $\alpha^k$ is a parameter tuned with a line search method. In
the implementation in Eigen~\cite{eigenweb}, $\lambda \ I_d$ is scaled with respect to the diagonal elements of
$(dG^{t, k-1}d\phi^{t, k-1})^T dG^{t, k-1}d\phi^{t, k-1}$. All the tolerances
for convergence
are set to machine precision, and the maximum number of residual evaluations is
set to $7$ unless explicitly stated differently in the numerical results section~\ref{sec:results}.

\begin{rmk}[Least-squares Petrov-Galerkin]
    \label{rmk:LSPG}
    The method is called manifold LSPG because it refers to the LSPG method
    usually applied when the manifold is linear. It consists in multiplying the
    residual to the left with a different matrix $\Psi$ with respect to the
    linear embedding $\Phi$ of the reduced coordinates into the state space,
    \begin{equation}
        \Psi^{T}G_{h, \delta t}(\boldsymbol{\mu}, \Phi \mathbf{z})=0,
    \end{equation}
    where, $\Phi\in\mathbb{R}^{d\times r}$ is the basis of the linear reduced
    manifold contained in $X_h\subset\mathbb{R}^{d}$,
    $\mathbf{z}\in\mathbb{R}^{r}$, and $\Psi$ is to be defined: the left
    subspace $\Psi\in\mathbb{R}^{d\times r}$ is used to enforce the
    orthogonality of the non-linear residual to a left subspace
    $\mathcal{L}\subset\mathbb{R}^{d}$. Applying Newton's method because of the
    nonlinearities, the problem is translated into the iterations for
    $k=1,\dots,K$:
    \begin{align}
        \Psi^{T}dG\Phi\ \delta\mathbf{z}^k=-\Psi^{T}G^k, \\
        \mathbf{z}^{k} =\mathbf{z}^{k-1} + \alpha^k\ \delta\mathbf{z}^{k}.
    \end{align}
    The step length $\alpha^k$ is computed after a line search along the
        direction $p^k$. Usually $\Phi$ is chosen from a POD basis of the state
        variable $\mathbf{z}\in X_h$. For the left subspace, that imposes
        orthogonality constraints, different choices can be applied. In general,
        given the system
    \begin{equation}
        \label{eq:modelII}
        dG\Phi\ \delta\mathbf{z}^k=-G^k,
    \end{equation}
    the least squares solution is the one orthogonal to the range of $dG\Phi$.
    In the case of $dG$ symmetric positive definite we have that $\Psi\subset
    <dG\Phi>=<\Phi>$ i.e. $\Psi=\Phi$ and the method is called Galerkin
    projection, but in general if this is not true then the optimal left
    subspace remains $\Psi=dG\Phi$. For examples where the Galerkin projection
    is not optimal see \cite{carlberg2011efficient}, section 3.4
    \textit{Numerical comparison of left subspaces}. This is often the case for
    advection-dominated discretized systems of PDEs: in these cases LSPG is
    preferred to Galerkin projection.\\
\end{rmk}

\begin{rmk}[Manifold Galerkin]
    The manifold Galerkin method proposed in~\cite{lee2020model} assumes that
    the columns of the Jacobian of the decoder are good approximations of the
    state velocity space: if a spatial discretization is applied and the
    residual has the form
    $\mathcal{G}_h(\boldsymbol{\mu}, \mathbf{U}) =
    \dot{\mathbf{U}}-\mathbf{f}(\boldsymbol{\mu}, \mathbf{U})$, where
    $\mathbf{f}$ is a generic, possibly non-linear, vector field, then
    \begin{equation}
        \dot{\mathbf{z}} = \argmin_{\mathbf{v}\in\mathbb{R}^r}\lVert d\phi(\mathbf{z})\mathbf{v}-\mathbf{f}(\boldsymbol{\mu}, \phi(\mathbf{z}))\rVert^2,
    \end{equation}
    is solved for the latent state velocity, under the hypothesis that
    $d\phi(\mathbf{z})$ has full-rank and $d\phi$ is a good approximation of the
    full-order state velocity even if the autoencoder is trained only on the
    values of the state for different times and parameters, without considering
    its velocty. If $\phi$ is linear,  we obtain the Galerkin method presented in the
    remark~\ref{rmk:LSPG}, after having applied a temporal discretization scheme
    and multplied the resulting equation to the left with $d\phi(\mathbf{z})=\Phi$ as
    is the case for linear Galerkin projection: the discretize-then-project and
    project-then-discretize approaches are equivalent in this
    case~\cite{lee2020model}.

    The LSPG performs better as shown
    in~\cite{lee2020model}, even if they are asymptotically equivalent also in
    the case of a
    non-linear manifold, provided $\phi$ is twice differentiable. So we have
    chosen to employ only the NM-LSPG method and do not compare it with the
    non-linear manifold Galerkin (NM-G) method.
\end{rmk}

For the numerical tests we have performed, the numerical approximation of the
Jacobian of the residual $G_{h, \delta t}(\boldsymbol{\mu}, \phi(\mathbf{z}),
\{\phi(\mathbf{z}^s)\}_{s\in I_t})$ is accurate enough. So in the
implementation, at each iteration step the Jacobian of the residual with respect
to the latent variable, that is $dG^{t, k-1}d\phi^{t, k-1}$ of
equation~\ref{eq:levenbergMarquardt}, is approximated with finite differences.
The step size is taken sufficiently lower than the distance between consecutive
latent states.
\subsection{Reduced over-collocation method}
At the point of equation~\ref{eq:levenbergMarquardt}, the model still depends on
the number of degrees of freedom of the full-order model $d$, since at each time
step and optimization step the latent reduced variable
$\mathbf{z}\in\mathbb{R}^r$ is forwarded to the reconstructed state
$\mathbf{U}_h = \phi(\mathbf{z})\in\mathbb{R}^d$. A possible solution is
represented by the reduced over-collocation method~\cite{chen2021eim}, for which
the least squares problem~\ref{eq:PDEsystemDiscrete} is solved only on a limited
number of points $r<r_h << d$,
\begin{equation}
    \mathbf{z}^t = \argmin_{\mathbf{z}\in \mathbb{R}^r}\ \lVert P_{r_h}G_{h, \delta t}(\boldsymbol{\mu}, \phi(\mathbf{z}), \{\phi(\mathbf{z}^s)\}_{s\in I_t})\rVert^{2}_{\mathbb{R}^{r_h}}
    \label{eq:PDEsystemDiscreteHyperReduced1},
\end{equation}
where $P_{r_h}$ is the projection onto $r_h$ standard basis elements in
$\mathbb{R}^{d}$ associated to the over-collocation nodes or magic points and
selected as described later. Afterwards the Levenberg-Marquardt algorithm is
applied as described in the previous section, to solve the least squares problem
~\ref{eq:PDEsystemDiscreteHyperReduced1}.

At this point we make the assumption that the method used to discretize the
model has a local formulation so that each discrete differential operator can be
restricted to the nodes/magic points of the hyper-reduction and consequently,
the least-squares problem~\ref{eq:PDEsystemDiscreteHyperReduced1} reduces to
\begin{align}
    \mathbf{z}^t &= \argmin_{\mathbf{z}\in \mathbb{R}^r}\ \lVert P_{r_h}G_{h, \delta t}(\boldsymbol{\mu}, P_{r_h}(\phi(\mathbf{z})), \{P_{r_h}(\phi(\mathbf{z}^s))\}_{s\in I_t})\rVert^{2}_{\mathbb{R}^{r_h}},\\
    &=\argmin_{\mathbf{z}\in \mathbb{R}^r}\ \lVert \Tilde{G}_{h, \delta t}(\boldsymbol{\mu}, \Tilde{\phi}(\mathbf{z}), \{\Tilde{\phi}(\mathbf{z}^s)\}_{s\in I_t})\rVert^{2}_{\mathbb{R}^{r_h}}
    \label{eq:PDEsystemDiscreteHyperReduced2},
\end{align}
where the projected residual $\Tilde{G}$ is introduced and the compressed
decoder $\Tilde{\phi}$ is defined to substitute $P_{r_h}\circ\phi$ with another
embedding from the latent space to the hyper-reduced space in $\mathbb{R}^{r_h}$,
such that the whole structure of the decoder is reduced as described in
subsection~\ref{subsec:compressedDecoderTraining} to further decrease the
computational cost.

The Levenberg-Marquardt method is applied also to the hyper-reduced system
\begin{align}
    \left((d\Tilde{G}^{t, k-1}d\Tilde{\phi}^{t, k-1})^T d\Tilde{G}^{t, k-1}d\Tilde{\phi}^{t, k-1} + \lambda\ I_d\right) \delta\mathbf{z}^{t,k} = -(d\Tilde{G}^{t, k-1}d\Tilde{\phi}^{t, k-1})^T \Tilde{G}^{t, k},\label{eq:levenbergMarquardtHyperReduced}\\
    \mathbf{z}^{t,k} = \mathbf{z}^{t,k-1} + \alpha^k\ \delta\mathbf{z}^{t,k},
    \label{eq:PDEHyperReductionLM}
\end{align}
and as for the manifold LSPG method, the Jacobian matrix $d\Tilde{G}^{t,
k-1}d\Tilde{\phi}^{t, k-1}$ is numerically approximated at each optimization
step in the implementations.

\begin{rmk}[Submesh needed to define the hyper-reduced differential operators]
    To compute $\Tilde{G}_{h, \delta t}$ in the nodes/magic points of the
    reduced over-collocation method, some adjacent degrees of freedom are
    needed by the discrete differential operators involved. So, actually, at
    each time step not only the values of the state variables at the
    magic points are needed, but also at the adjacent degrees of freedom in the mesh with
    possible overlappings. We represent the restriction to this submesh of magic
    points and adjacent degrees of freedom with the projector
    $P^{s}_{r_h}\in\mathbb{R}^{s_h\times d}$, where $s_h$ is the number of
    degrees of freedom of the submesh. Equation~\ref{eq:PDEsystemDiscreteHyperReduced2} becomes
    \begin{equation}
        \label{eq:submesh residual}
        \mathbf{z}^t = \argmin_{\mathbf{z}\in \mathbb{R}^r}\ \lVert P_{r_h}G_{h, \delta t}(\boldsymbol{\mu}, P^{s}_{r_h}(\phi(\mathbf{z})), \{P^{s}_{r_h}(\phi(\mathbf{z}^s))\}_{s\in I_t})\rVert^{2}_{\mathbb{R}^{r_h}}.
    \end{equation}
    The stencil around each magic point to
    consider
    depends on the type of numerical scheme. Since we are using the finite
    volume method we have to consider the degrees of freedom of the adjacent
    cells. For example, for cartesian grids, the schemes chosen for
    the 2d non-linear conservation law have a stencil of $1$ layer of adjacent
    cells, $4$ additional nodes in total for a cell of the interior of the mesh. The 2d shallow water equations case requires a stencil
    of $2$ layers instead, $12$ additional nodes in total for a cell of the interior of the mesh. The two cases are
    shown in Figure~\ref{fig:mp stencils}.

    \begin{figure}[h]
        \centering
        \includegraphics[width=0.4\textwidth]{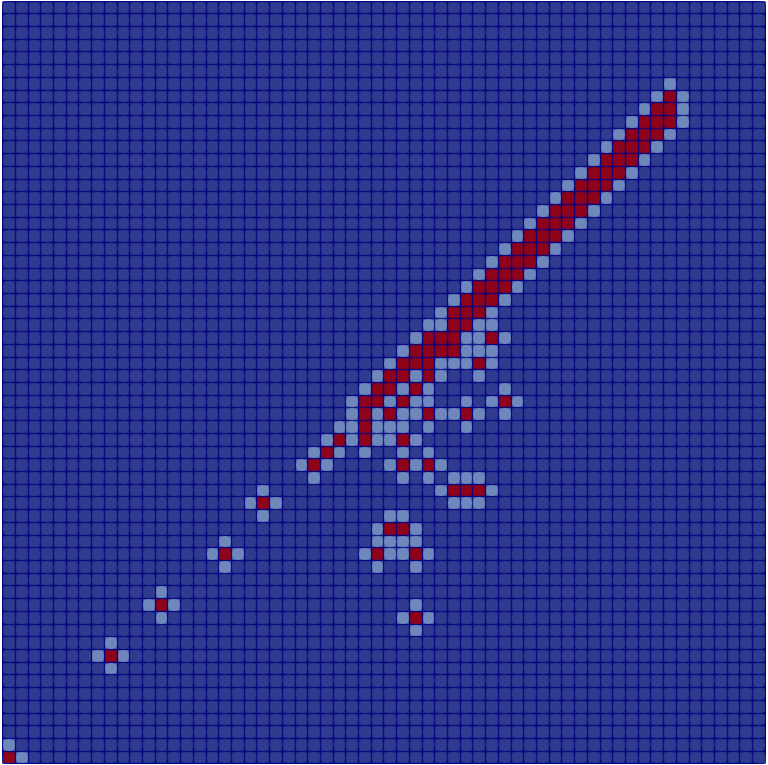}
        \includegraphics[width=0.4\textwidth]{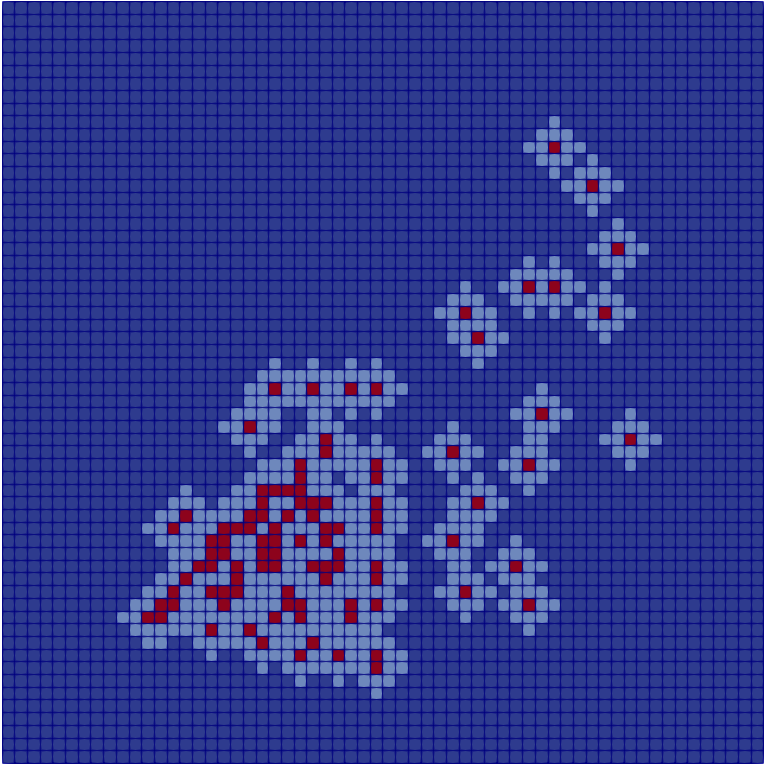}
        \caption{Left: 100 magic points of the 2d non-linear conservation law test case. Right: 100 magic points of the 2d shallow water test case. The magic points are represented in red, the stencils of the cells and associated degrees of freedom needed for the evaluation of the discrete differential operators are in light-blue. The discarded nodes in the evolution of the dynamics with NM-LSPG-ROC are in blue. The stencil is made by $1$ layer of cells in the NCL case and $2$ layers in the SW case.}
        \label{fig:mp stencils}
    \end{figure}
\end{rmk}

\subsection{Over-collocation nodes selection}
The nodes/magic points of the over-collocation hyper-reduction method should be
defined such that
\begin{equation}
    P_{r_h} = \argmin_{P^T P\in\mathcal{S}}\ \max_{(\mathbf{z}^t, \{\mathbf{z}^s\}_{s\in I_t})\in\mathcal{T}}\ \lVert G_{h, \delta t}(\boldsymbol{\mu}, (\phi(\mathbf{z})), \{(\phi(\mathbf{z}^s))\}_{s\in I_t}) - P^TPG_{h, \delta t}(\boldsymbol{\mu}, P(\phi(\mathbf{z}^t)), \{P(\phi(\mathbf{z}^s))\}_{s\in I_t})\rVert^{2}_{\mathbb{R}^{r_h}}
    \label{eq:minimizationHyperReductionProjector}
\end{equation}
where $\mathcal{S}=\{P\in\mathbb{R}^{r_h\times
d}|P=(\mathbf{e}_{i_1}|\dots|\mathbf{e}_{i_{r_h}})^T\}$ is the space of
projectors onto $r_h$ coordinates associated to the standard basis
$\{\mathbf{e}_i\}_{i\in\{1,\dots,d\}}$ of $\mathbb{R}^d$ and $\mathcal{T}$ is
the space of discrete solution trajectories varying with respect to
$\boldsymbol{\mu}$ and the intermediate optimization steps
\begin{align*}
    \mathcal{T}&=\{(\boldsymbol{\mu}, t, k, \mathbf{z}^{t, k}, \{\mathbf{z}^{s, k}\}_{s\in I_t})\in\mathcal{P}\times V_{h, \boldsymbol{\mu}}\times V_{h, \boldsymbol{\mu},
    t}\times\mathbb{R}^d\times\mathbb{R}^{d\times |I_t|}|,\\
&\left((dG^{t, k-1}d\phi^{t, k-1})^T dG^{t, k-1}d\phi^{t, k-1} + \lambda\ I_d\right) \delta\mathbf{z}^{t,k} = -(dG^{t, k-1}d\phi^{t, k-1})^T G^{t, k}\},
\end{align*}
where $V_{h, \boldsymbol{\mu}}$ is the discrete space of time instants, possibly
depending on $h$ and the paremter $\boldsymbol{\mu}$, and $V_{h,
\boldsymbol{\mu}, t}$ is the discrete space of optimization steps at time $t$.
Essentially we want that the nodes/magic points approximate the residuals among
all time steps, optimization steps and parameter instances.

There are many possible algorithms to solve
Equation~\ref{eq:minimizationHyperReductionProjector} for $P_{r_h}$. Usually they
are not optimal and compromise between computational cost and accuracy,
depending on the problem at hand. Among others, these algorithms are part of the
hyper-reduction methods such as empirical interpolation
method~\cite{barrault2004empirical}, discrete empirical interpolation
method~\cite{chaturantabut2010nonlinear}, Gauss-Newton tensor approximation (GNAT)~\cite{carlberg2013gnat},
space-time GNAT~\cite{choi2019space} and solution-based non-linear subspace GNAT~\cite{choi2020sns}(SNS-GNAT).

In particular, if $\mathcal{G}_h(\boldsymbol{\mu}, \mathbf{U}) =
\dot{\mathbf{U}}-\mathbf{f}(\boldsymbol{\mu}, \mathbf{U})$, then, following some
considerations that justify SNS-GNAT~\cite{choi2020sns}, the training modes
employed to find the nodes/magic points of the over-collocation method are
represented by the state snapshots instead of the residual fields. We could in
principle use the reduced fields $\phi(\mathbf{z})$ but in practice, for the
test case we considered, the full-order state snapshots were enough, without
even saving the intermediate optimization states.

The procedure is applied at the same time for all the components of the state
field $\mathbf{U}\in X_h$ and follows a greedy approach. We remark that
the
spatial discretization must not vary, so that the degrees of freedom correspond
to the same spatial and physical quantity over time and for every parameter
instance. Some approaches tackle also geometry deformations, but keeping the
same number of degrees of freedom in a reference
system~\cite{stabile2020efficient}.

The algorithm~\ref{algo:nodesRoc} is an adaptation of GNAT from Algorithm 3
in~\cite{Carlberg_2013} to the simpler case in which the Jacobian matrix is not
considered in the hyper-reduction (since in the LM method the Jacobian matrix is
approximated with finite differences from the residual, see
equation~\ref{eq:PDEHyperReductionLM}). Also with respect to Algorithm 3
in~\cite{Carlberg_2013}, the new node/magic point at line 16 in
Algorithm~\ref{algo:nodesRoc} is found without computing also the reconstruction
error of the degrees of freedom associated to its stencil.

\IncMargin{1em}
\begin{algorithm}[htb]
    \DontPrintSemicolon \caption{Greedy nodes/magic points evaluation for ROC
        method.}\label{alg:roc}
        \SetKwInOut{Input}{input}\SetKwInOut{Output}{output}

    \Indm \Input{training state fields
    $\mathcal{U}_{\text{train}}:=\{\mathbf{U}_{\boldsymbol{\mu}, t}\}_{\mu\in
    \mathcal{P}_{\text{train}},\ t\in V_{\boldsymbol{\mu}, h}}$, such that
    $|\mathcal{U}_{\text{train}}|=n_{\text{train}}$,\\
    $r_{\text{init}}$ nodes/magic points at the boundaries and forcefully
    added,\\
     number of nodes/magic points $r_h$,\\
     $N_{\text{modes}}$ number of training modes extracted to compute
     nodes/magic points.} \Output{ $r_h$ nodes/magic points used to define
     $P_{r_h}$.}

    \Indp \BlankLine Extract $N_{\text{modes}}$ modes
    $\{\phi^{i}_{\mathbf{U}}\}_{i\in\{ 1,\dots, N_{\text{modes}}\}}$ from the
    SVD of $\mathcal{U}_{\text{train}}$.\\
    Compute the additional number of nodes to sample $n_a =
    r_h-r_{\text{init}}$.\\
    Initialize counter for the number of working basis vectors used: $n_b =
    0$.\\
    Set the number of greedy iterations to perform: $n_{it} = \min{n_c, n_a}$.\\
    Compute the minimum number of working basis vectors per iteration:
    $n_{\text{ci}, \text{min}}=\text{floor}(n_c/n_it)$.\\
    Compute the minimum number of sample nodes to add per iteration:
    $n_{\text{ai}, \text{min}}=\text{floor}(n_a/n_c)$.\\
    \For(greedy iteration loop){$i=1,\dots,n_{it}$}{Compute the number of workin
    basis vectors for this iteration: $n_{ci} = n_{\text{ci}, \text{min}}$; If
    $i\leq n_c,\ \text{mod}\ n_{it}$ then $n_{ci}=n_{ci}+1$.\\
    Compute the number of samples nodes to add during this iteration: $n_{ai} =
    n_{\text{ci}, \text{min}}$; If $i\leq n_a,\ \text{mod}\ n_{c}$ then
    $n_{ai}=n_{ai}+1$.\\
    \eIf{$i=1$}{$\left[\mathbf{U}^{1}\cdots \mathbf{U}^{{n}_{ci}}\right] =
    \left[\phi^{1}_{\mathbf{U}}\cdots \phi^{n_{ci}}_{\mathbf{U}}\right]$,\\
    build the intermediate projection $P$ from the $r_{\text{init}}$ initial
    nodes/magic points.}{\For(basis vector
    loop){$q=1,\dots,n_{ci}$}{$\mathbf{U}^{q} =
    \phi^{n_{q}+q}_{\mathbf{U}}-\left[\phi^{1}_{\mathbf{U}}\cdots
    \phi^{n_b}_{\mathbf{U}}\right]\alpha$, with
    $\alpha=\argmin_{\gamma\in\mathbb{R}^{n_b}}\lVert\left[P\phi^{1}_{\mathbf{U}}\cdots
    P\phi^{n_b}_{\mathbf{U}}\right]\gamma-P\phi^{n_b
    +q}_{\mathbf{U}}\rVert_{2}$}} \For(sample node
    loop){$j=1,\dots,n_{ai}$}{Find the node/magic point $n$ that maximize the
    state reconstruction error among the nodes not selected yet: $n=\argmax
    \sum_{q=1}^{n_{ci}}(U^{q})^2$.\\
    Update $P$ with the newly found node/magic point.} $n_b = n_b + n_{ci}$}
\label{algo:nodesRoc}
\end{algorithm}
\DecMargin{1em}

\begin{rmk}[Comparison beetween GNAT and reduced over-collocation]
    We must say that the GNAT method, employed for the implementation of NM-LSPG
    with shallow masked autoencoders in~\cite{kim2020fast}, is a generalization of the reduced
    over-collocation method. In some cases though, they may perform similarly.
    For $P\in\mathcal{S}$, let us define
    \begin{equation}
        \mathbf{r}=r_{t, \boldsymbol{\mu}, k} = G_{h, \delta
    t}(\boldsymbol{\mu}, P(\phi(\mathbf{z}^{t, k})), \{P(\phi(\mathbf{z}^{s, k-1}))\}_{s\in
    I}),\ \forall (t, k, \boldsymbol{\mu})\in V_{h, \boldsymbol{\mu}}\times
    V_{h, \boldsymbol{\mu}, t}\times \mathcal{P},
    \end{equation}
    and the GNAT projection operator $\mathbb{P}=\Phi (P\Phi)^{\dagger}P$,
    where $\Phi$ is the matrix where the columns correspond to a chosen basis
    (it could be the FOM residual snapshots, ROM residual snapshots, ROM
    jacobian and residual snapshots, see~\cite{carlberg2013gnat}). In particular,
    if $\Phi=P_{r_h}^{T}=P^{T}$ we have that
    \begin{equation}
        \mathbb{P}=\Phi (P_{r_h}\Phi)^{\dagger}P_{r_h}=\Phi (P_{r_h}\Phi)^{-1}P_{r_h}=P_{r_h}^T (P_{r_h}P_{r_h}^T)^{-1}P_{r_h}=P_{r_h}^{T}P_{r_h},
    \end{equation}
    that is the reduced
    over-collocation projection in the chosen nodes/magic points and extended to
    $0$ in the remaining degrees of freedom. In this sense, the GNAT method
    includes the reduced over-collocation one.

    However, for the class of problems we are considering, the GNAT method suffers from the slow decaying KnW. We have the inequalities
    \begin{align*}
        d_{n}^2(\{r_{t, \boldsymbol{\mu}, k}\}_{t, \boldsymbol{\mu}, k})&= \inf_{\text{dim}V_n=n}\max_{t, \boldsymbol{\mu}, k} \lVert \mathbf{r}-V_n V_n^T \mathbf{r}\rVert^2_2 \leq \inf_{\text{dim}V_n=n}\max_{t, \boldsymbol{\mu}, k} \lVert \mathbf{r}-V_n V_n^T \mathbf{r}\rVert^2_2 + \lVert V_n V_n^T \mathbf{r} - \mathbb{P}\mathbf{r}\rVert^2_2\\
        &=\inf_{\text{dim}V_n=n}\max_{t, \boldsymbol{\mu}, k} \lVert \mathbf{r}- \mathbb{P}\mathbf{r}\rVert^2_2 = \inf_{\text{dim}V_n=n}\max_{t, \boldsymbol{\mu}, k} \lVert (I- \mathbb{P})(I-V_nV_n^T)\mathbf{r}\rVert^2_2\\
        &\leq \inf_{\text{dim}V_n=n}\max_{t, \boldsymbol{\mu}, k} \lVert (I- \mathbb{P})\rVert_{2}^{2}\lVert(I-V_nV_n^T)\mathbf{r}\rVert^2_2 = \inf_{\text{dim}V_n=n}\max_{t, \boldsymbol{\mu}, k} \lVert \mathbb{P}\rVert_{2}^{2}\lVert(I-V_nV_n^T)\mathbf{r}\rVert^2_2,
    \end{align*}
    where $(t, k, \boldsymbol{\mu})\in V_{h, \boldsymbol{\mu}}\times
    V_{h, \boldsymbol{\mu}, t}\times \mathcal{P}$ whenever the maximum is taken.
    The term $\mathbf{r}- \mathbb{P}\mathbf{r}$ is the GNAT approximation error. The rightmost term is the usual bound on the hyper-reduction
    error~\cite{chaturantabut2010nonlinear}, where it was used the fact that $\mathbb{P}=I-\mathbb{P}$.
    The second equality is valid because $\mathbf{r}-V_n V_n^T \mathbf{r}$ and
    $V_n V_n^T \mathbf{r} - \mathbb{P}\mathbf{r}$ are orthogonal. The third
    equality is obtained from the relations
    \begin{align*}
        &\mathbf{r}=(\mathbf{r}-\mathbf{r}_{*})+\mathbf{r}_{*}=\mathbf{w}+\mathbf{r}_{*},\qquad \text{with}\ \mathbf{r}_{*}:=V_n V_n^T \mathbf{r},\\
        \Rightarrow \quad&\mathbf{r}-\mathbb{P}\mathbf{r}=\mathbf{w} + \mathbf{r}_{*} - \mathbb{P}\mathbf{w}- \mathbb{P}\mathbf{r}_{*} = \mathbf{w}- \mathbb{P}\mathbf{w},
    \end{align*}
    where the last equality follows from $\mathbb{P}\mathbf{r}_{*}=\mathbf{r}_{*}$.
    Here we have supposed the nodes/magic points to be independent from $V_n$.
    So in the case of slow decaying Kolmogorov n-width, minimizing the
    hyper-reduced residual $\mathbb{P}\mathbf{r}\in<V_n>\subset\mathbb{R}^d$ is less efficient due to the slow
    convergence in $n$ of the best approximation error $\lVert\mathbf{r}-V_n
    V_n^T\mathbf{r}\rVert^{2}_2$. This is one of the reasons why we employed ROC
    for the SWE test case; for the NCL test case GNAT and ROC performed similarly.

\end{rmk}
\subsection{Compressed decoder teacher-student training}
\label{subsec: compressed decoder training}
In order to make the whole methodology independent from the number of degrees of
freedom, the decoder has to be substituted with a map
$\Tilde{\phi}:\mathbb{R}^{R}\rightarrow P_{r_h}(X_h)\subset\mathbb{R}^{r_{h}}$ from the latent space
to the space of discrete full-order solutions evaluated only at the submesh
containing the magic points and the needed adjacent degrees of freedom. As architecture, we choose a feedforward neural
network (FNN) with one hidden layer, but actually the only requirement is that the
computational cost is low enough such that not only a theoretical dimension
reduction is achieved but also a speed-up is reached.

\begin{figure}[h]
    \centering
    \includegraphics[width=0.8\textwidth]{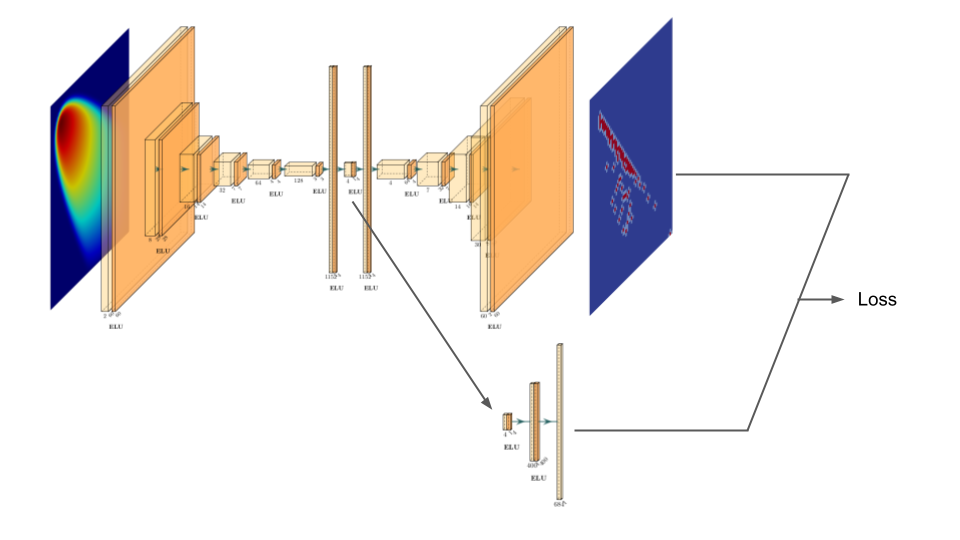}
    \caption{Teacher-student training of the compressed decoder for the 2d
    non-linear conservation law test case.
    The magic points, which the snapshots are restricted to, are shown in red over
    the domain.}
    \label{fig:teacher student net}
\end{figure}

In the literature, the procedure for the training
of the compressed decoder $\Tilde{\phi}$ from the decoder $\phi$ is called
teacher-student training~\cite{gou2021knowledge}. In principle, the compressed decoder can be composed
from layers inherited by the original decoder, such that the learning process
involves only the final new additional layers. In our case, we preferred to
train the compressed decoder anew: the latent projections of the training snapshots with
the encoder $\psi(\mathcal{U}_{\text{train}})=\{z_{i}\}_{i=1}^{N_{\text{train}}}$ are the inputs and the restriction of
the snapshots to the submesh $P^s_{r_h}(\mathcal{U}_{\text{train}})=\{\Tilde{U}_{i}\}_{i=1}^{N_{\text{train}}}$ are the
targets, see Equation~\ref{eq:submesh residual}. A schematic representation of the teacher-student training is
represented in figure~\ref{fig:teacher student net}. Moreover, to speed-up the
offline tage, we use the training FOM snapshots restricted to the magic points
as training outputs for the teacher-student training, while usually the
reconstructed snapshots from the CAE -that would additionaly need to be
computed- are employed.

\label{subsec:compressedDecoderTraining}

Again, for the training,
we use a relative mean square loss with an additional regularizing term
\begin{equation}
\label{eq: decoder compress}
    \mathcal{L}(\{z_{i}\}_{i=1}^{b};\Tilde{\Theta}) = \frac{1}{b}\sum_{i=1}^{b}\frac{\lVert \Tilde{U}_{i} - \Tilde{\phi}(z_{i})\rVert^{2}_{2}}{\lVert \Tilde{U}_{i}\rVert^{2}_{2}} + \lambda_{1} \lVert \Tilde{\Theta} \rVert^{2}_{2},
\end{equation}
where $\Tilde{\Theta}$ are the weights of the compressed decoder.

\begin{rmk}[Jacobian evaluation in Levenberg-Marquardt algorithm]
    The main reason why finite differences approximations of Jacobians are
    implemented in the NM-LSPG case, as explained at the end of
    subsection~\ref{subsec: NM-LSPG}, is that the computational cost of
    evaluating the Jacobian of the full decoder is too high. In principle
    Jacobian evaluations of the compressed decoder are cheaper and could be
    employed, instead of relying again on
    finite differences approximations.
\end{rmk}

\begin{rmk}[Shallow masked autoencoders]
    We are motivated to write this article to extend the results
    in~\cite{kim2020fast} to a generic architecture composed by neural networks.
    They performed the hyper-reduction of the non-linear manifold
    method~\cite{lee2020model} with a shallow masked autoencoder, so that
    correctly masking the weights matrices of the decoder, its outputs
    correspond only to the submesh needed by the GNAT method, thus eliminating
    the dependence on the FOM's degrees of freedom. We want to reproduce, in
    some sense, this approach for
    an arbitrary autoencoder architecture, in this case a CAE, in order to takle with the latest architectures developed in
    the literature the
    problem of solution manifold approximability: we think this is a major
    concer when trying to apply non-linear MOR to real applications. In fact, as
    will be clear in the numerical results section~\ref{sec:results}, the reconstruction error of the
    autoencoder bounds from below the prediction error of our newly developed ROMs.
\end{rmk}

It can be seen that the new model order reduction is composed of
two distinct procedures to achieve the independency on the number of degrees of
freedom: first the residual from NM-LSPG in Equation~\ref{eq:PDEsystemDiscrete}
is hyper-reduced with ROC in Equation~\ref{eq:PDEsystemDiscreteHyperReduced1}
and secondly the CAE's decoder is compressed with teacher-student training. In
principle, we could substitute the use of the compressed decoder with the
restriction of the final layer of the CAE's decoder into the magic points, while
keeping the hyper-reduction with ROC of the residual. In this case, the whole
methodology would still be dependent on the total number of degrees
of freedom, but in practice a CAE's decoder forward is relatively cheap compared
to the evaluation of the full residual. So, the hyper-reduction performed with
ROC or GNAT only at the equations/residuals level, is already beneficial to
reduce the computational cost. We will compare this variant of the NM-LSPG-ROC
method with the one that employs the compressed decoder, also to verify the
consistency of the teacher-student training that is omitted in the first case.

In the numerical results section~\ref{sec:results} we will adopt the acronym
NM-LSPG-ROC-TS or NM-LSPG-GNAT-TS for the method that employs the compressed
decoder and NM-LSPG-ROC or NM-LSPG-GNAT for the method that performs the
hyper-reduction only at the equations/residuals level.

\section{Numerical results}
\label{sec:results}
We test the new methodology on two benchmarks with a relatively slow KnW: the first
model is governed by a
non-linear conservation law~\ref{subs:nonlinear conservation law} the second by
the shallow water equations~\ref{subs:shallow water}. Both are parametric, non-linear and
time-dependent, and the only other (non-temporal) parameter is a multiplicative constant of the initial
conidition. The mesh employed is the same: a $60\times 60$ structured orthogonal
grid.

All the CFD simulations are obtained by the use of an in-house open
source library ITHACA-FV (In real Time Highly Advanced Computational Applications
for Finite Volumes)~\cite{ITHACA-FV}, developed in a finite volume environment
based on the open-source library
OpenFOAM~\cite{OF}. Regarding the implementation of the convolutional
autoencoders and compressed decoders (CAE) we used LibTorch, PyTorch C++ frontend,
while for the training of the long-short term memory network (LSTM) we used
PyTorch~\cite{NEURIPS2019_9015}. All the CFD simulations were performed on a
Intel(R) Core(TM) i7-8750H CPU with 2.20GHz and all the neural networks
trainings on a GeForce GTX 1060 GPU. Further
reductions in
the computational costs could be achieved exploiting the parallel implementation
of the training procedures in PyTorch. The details of the architectures of the neural networks that will be employed are reported in
the Appendix~\ref{sec:NN}.

The following notations are introduced: $N^{\mu}_{\text{train}},\
    N^{\mu}_{\text{test}}$ are the numbers of train and test parameters,
    respectively; $N^{t}_{\text{train}},\ N^{t}_{\text{test}}$ are the number of
    time instances associated to the train and test parameters, respectively.
    The total number of training and test snapshots is thus
    $N_{\text{train}}=N^{\mu}_{\text{train}}\cdot N^{t}_{\text{train}}$, and
    $N_{\text{test}}=N^{\mu}_{\text{test}}\cdot N^{t}_{\text{test}}$,
    respectively.

The accuracy of the reduced-order models devised is measured with the mean
realtive $L^{2}$-error and the maximum relative $L^{2}$-error, where the mean and
max are taken with respect to the time scale: since the test cases depend on a
non-temporal parameter, for each instance of these parameters a time-series
corresponding to the discrete dynamics is associated; the mean and maximum are
evaluated w.r.t the elments of these time-series. Let $\{u^{t_i}_{\mu}\}_{i=1,\dots
    N^{t}}$ and $\{U^{t_i}_{\mu}\}_{i=1,\dots
    N^{t}}$ be the predicted and true time-series $N^{t}$ elements long, associated to the train or test
parameter $\mu$, the mean relative $L^2$-errors and maximum relative $L^2$-errors
are then defined as
\begin{equation}
    \epsilon_{mean}(u_{\mu}, U_{\mu}) = \frac{1}{N^{t}}\sum_{i=1}^{N^{t}} \frac{\lVert u^{t_i}_{\mu}-U^{t_i}_{\mu}\rVert_{L^{2}}}{\lVert U^{t_i}_{\mu}\rVert_{L^{2}}},\qquad \epsilon_{max}(u_{\mu}, U_{\mu}) = \max_{i=1,\dots,N^{t}} \frac{\lVert u^{t_i}_{\mu}-U^{t_i}_{\mu}\rVert_{L^{2}}}{\lVert U^{t_i}_{\mu}\rVert_{L^{2}}}.
\end{equation}

\begin{rmk}[Levenberg-Marquardt parameters]
    \label{rmk: lm parameters}
    Regarding the Levenberg-Marqurdt non-linear optimization algorithm, we remark
    that we approximate the Jacobians with forward finite differences, and the
    optimization process, for each time step, is stopped when the
    maximum number of residual evaluations is reached. This number is set to $7$, including the evaluations
    related to the Jacobian computations. When $7$ residual
    evaluations are not enough for the method to converge, it is explicitly
    reported.
\end{rmk}

\subsection{Non-linear conservation law (NCL)}
\label{subs:nonlinear conservation law}
We test our procedure for non-linear model order reduction on a 2d non-linear
conservation law model (NCL). Two main reasons are behind this choice: the slow
Kolmogorov
n-width decay of the continuous solution manifold, and
the possibility to compare our results with a similar test case realized with an implementation of non-linear
manifold based on shallow masked autoencoders and GNAT~\cite{kim2020fast}.

The parametrization affects the initial velocity as a scalar multiplicative
constant $\mu\in[0.8, 2]$:
\begin{equation}
    \label{eq:burgers 2d}
    \begin{cases}
        \partial_{t} \mathbf{u} + \frac{1}{2}\nabla\cdot\left(\mathbf{u}\otimes\mathbf{u}\right) = \nu\Delta \mathbf{u} & (\mathbf{x}, t)\in[0, 1]^{2}\times [0, 2],        \\
        \mathbf{u}(\mathbf{x}, 0) = 0.8 \cdot \mu \cdot \sin(2\pi x)\sin(2\pi y)\chi_{[0, 0.5]^2}                       & \mathbf{x}\in[0, 1]^{2},                          \\
        \mathbf{u}(\mathbf{x}, t) = 0                                                                                   & (\mathbf{x}, t)\in\partial[0, 1]^{2}\times[0, 2],
    \end{cases}
\end{equation}
where the viscosity $\nu=0.0001$. We will collect $N^{\mu}_{\text{train}}=12$ equispaced training
parameters from the range $\mu\in[0.8, 2]$ and $N^{\mu}_{\text{test}}=16$ equispaced test parameters
from the range $\mu\in[0.6, 2.2]$. The first two and the last two parameters
will account for the extrapolation error. The time step is equal to
$\Delta t = 1\mathrm{e}-3$ seconds, but the training snapshots are collected every $4$
timesteps and the test snapshots every $20$, thus $N^{t}_{\text{train}}=501$,
and $N^{t}_{\text{test}}=101$. In the predictive online phase, the
dynamics will be evolved with the same timestep $\Delta t = 1\mathrm{e}-3$. For easyness of representation, the train
parameters are labelled from $1$ to $12$, and the
test parameters are labelled from $1$ to $16$.

\begin{figure}[h]
    \includegraphics[width=0.24\textwidth]{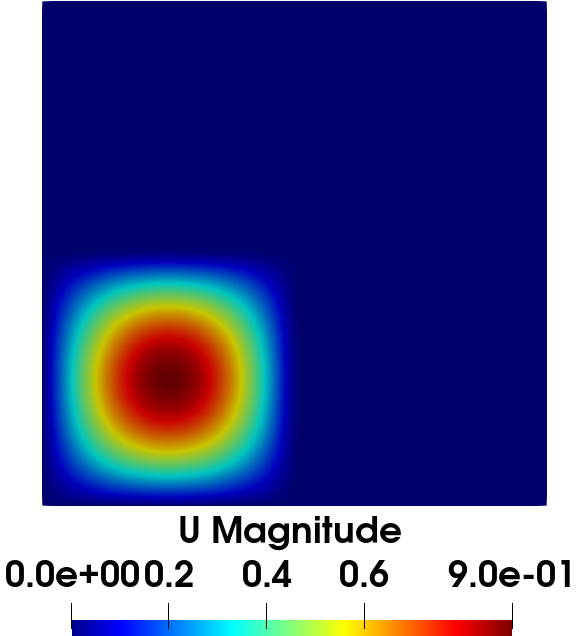}
    \includegraphics[width=0.24\textwidth]{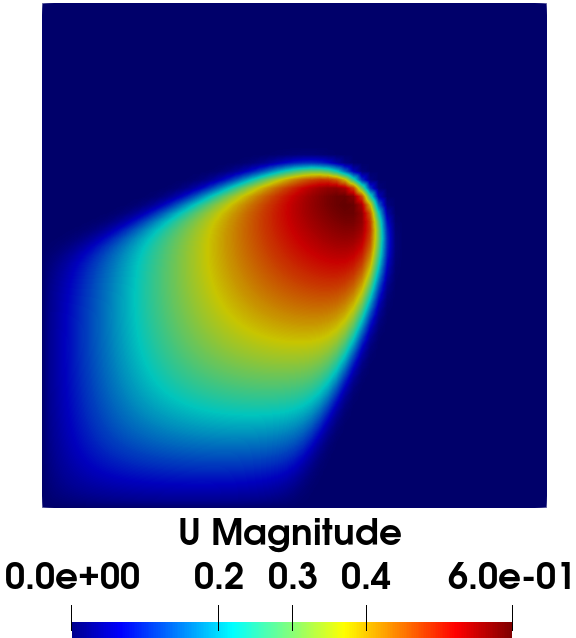}
    \includegraphics[width=0.24\textwidth]{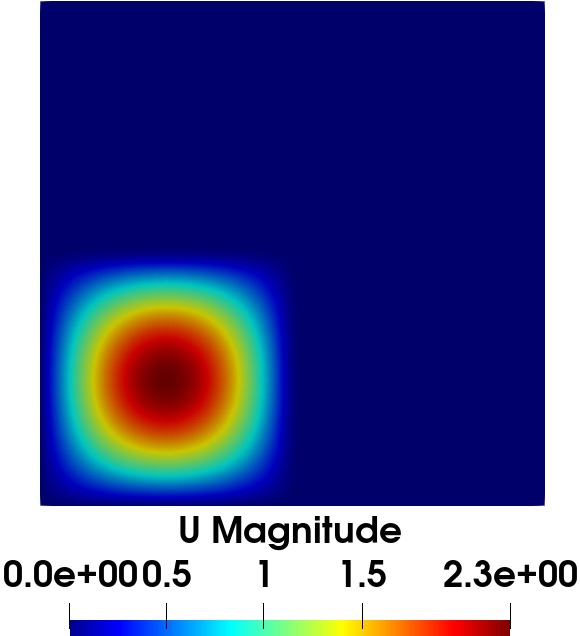}
    \includegraphics[width=0.24\textwidth]{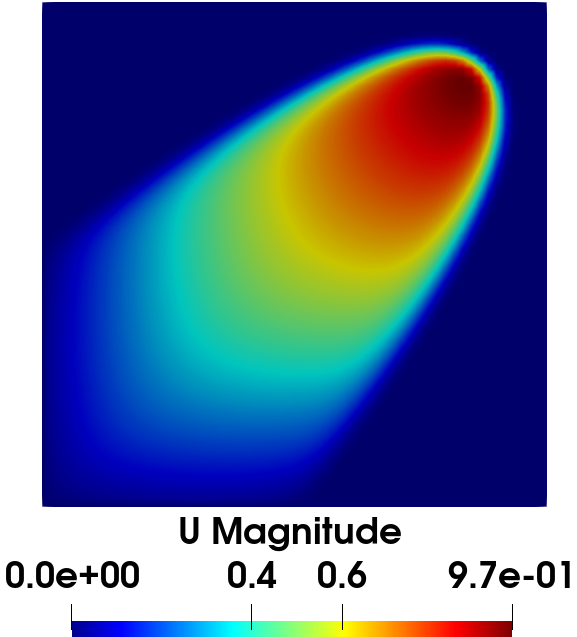}
    \caption{From left to right: FOM solution of equation~\ref{eq:burgers 2d} at
        $(t, \mu)\in\{(0, 0.8), (2, 0.8), (0, 2),(2, 2)\}$.}
    \label{fig:qualitative snap ncl}
\end{figure}

To have a qualitative view on the range of the solution manifold, we report
the initial and final time snapshots for the extremal training parameters of the range
$\mu\in[0.8, 2]$, in Figure~\ref{fig:qualitative snap ncl}.

In this test case the GNAT method performed slightly better than the ROC method for
hyper-reduction so we employed the former to obtain the results shown.

\subsubsection{Full-order model}
We solve the 2d non-linear conservation law for different values of the parameter
$\mu$ with OpenFoam~\cite{OF} open-source software for CFD. We employ the finite volumes
method (FVM) in a structured orthogonal grid of $60\times 60$ cells. If we represent
with $M$ the mass matrix, with $D$ the diffusive matrix term, and with $C(U^{t-1})$ the
advection matrix, then, at every time
instant $t$, the discrete equation
\begin{equation}
    \frac{M}{\Delta t}U^t+C(U^{t-1})U^t - \nu DU^t  = \frac{M}{\Delta t}U^{t-1},
\end{equation}
is solved for the state $U^t$ with a semi-implicit Euler method. The time step
is $1\mathrm{e}-3$, the initial and final time instants are $0$ and $2$ seconds.
The linear system is solved with the iterative
method BiCGStab preconditioned with DILU, until a tolerance of
$1\mathrm{e}{-17}$ on the
FVM residual is reached.

The stencil of the numerical scheme at
each cell involves the adjacent cells that share an interface ($4$ for an
interior cell, $3$ for a boundary cell and $2$ for a corner cell): the value of the state at the
interfaces is obtained with the bounded upwind method for the advection term and
the surface normal gradient is obtained with central finite differences of two
adjacent cell centers. So, in order to implement the reduced over-collocation
method, for each node/magic point we have to consider an additional number of
maximum $4$ cells, that is $8$ degrees of freedom to keep track of during the
evolution of the latent dynamics; of course in practice they may overlap reducing
the computational cost furtherly.

The residual of the NM-LPSG methods is evaluated with the same numerical scheme
of the FOM. In the SWE test case the FOM and the ROMs employ different numerical
schemes~\ref{subsec: SWE}.

\subsubsection{Manifold learning}
As first step of the procedure the discrete solution manifold is learned through
the training of a convolutional autoencoder (CAE) whose specific architecture is
reported in Table~\ref{tab: dec ncl}. The CAE is trained with the ADAM~\cite{kingma2014adam}
stochastic optimization algoritm for $2000$ epochs, halving the learning rate by
a factor of $2$ if after $200$ epochs the loss does not decrease. The initial
learning rate is $1\mathrm{e}{-3}$, its lower bound is $1\mathrm{e}{-6}$. The
number of training snapshots is $N_{\text{train}}=12\times 501=6012$, the batch size $20$. It could be furtherly refined in order the increase the efficiency of the
whole procedure.

We choose as latent dimension $4$, $2$ dimensions greater than the number of
parameters (the scalar multiplying the initial condition and time). We don't
perform a convergence study of the accuracy with respect to the latent
dimension since our focus is on the implementation of the NM-LSPG-ROC and
NM-LSPG-ROC-TS model-order reduction methods: we are satisfied as long as the
accuracy is relatively high, while the reduced dimension corresponds to an inaccurate
linear approximating manifold spanned by the same number of POD modes.

\begin{figure}[htp!]
    \centering
    \includegraphics[width=0.8\textwidth]{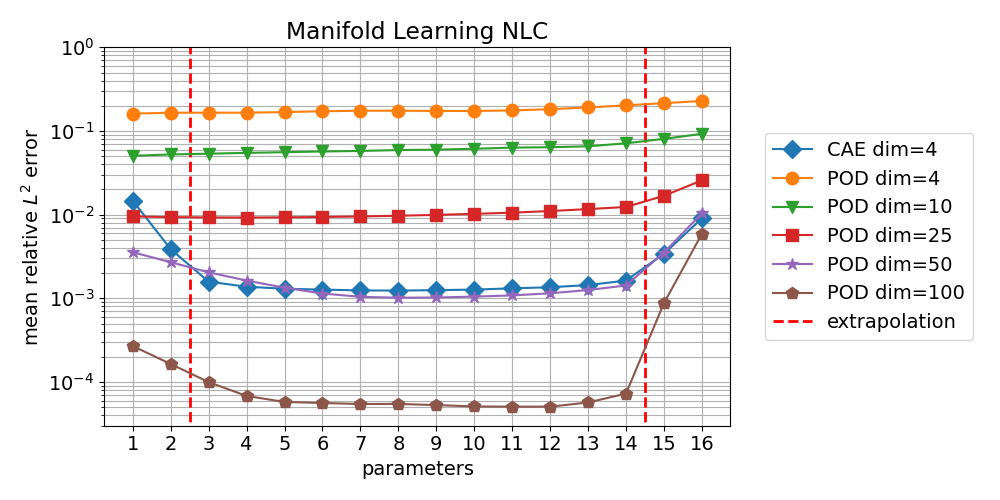}
    \caption{Comparison between the CAE's and POD's projection errors of the
    discrete solution manifold represented by the
    $N_{\text{train}}=N^{\mu}_{\text{train}}\cdot N^{t}_{\text{train}}=12\cdot 501$ training
    snapshots. The error is evaluated on the $16$ test parameters, each
    associated with a time series of $N^{t}_{\text{test}}=101$, for a total of
    $N_{\text{test}}=N^{\mu}_{\text{test}}\cdot N^{t}_{\text{test}}=16\cdot 101
    = 1616$ test
    snapshots. The mean is performed over the time scale.}
    \label{fig: ml_ncl}
\end{figure}

In
Figure~\ref{fig: ml_ncl} is shown the reconstruction error of the CAE and its
decay with respect to the number of POD modes chosen $[4, 10, 25, 50, 100]$. To
reach the same accuracy of the CAE with latent dimension $4$, around $50$ POD
modes are needed. In order to state that the slow KnW decay problem is overcome
by the CAE, the asymptotic convergence of the reconstruction error w.r.t. the
latent dimension should be studied as was done for similar problems
in~\cite{lee2020model,kim2020fast}. Instead, we will empirically prove that we can devise
an hyper-reduced ROM with latent dimension $4$ and accuracy lower than the $2\%$ for the
mean relative $L2$-error, a task that would be impossible for a POD based ROM
with the same reduced dimension, since the reconstruction error is near $20\%$
for all the test parameters.

Without imposing any additional inductive bias a part from the regularization
term in the loss from Equation~\ref{eq:loss cae} and the positiveness of the velocity
components, the latent trajectories reported in Figure~\ref{fig:
    latent_true_ncl} for the odd parameters of the test set, are qualitatively smooth. An
important detail to observe is that the initial conditions are well separated
one from another in the latent space, see Remark~\ref{rmk:initial condition}, and that the dynamics is non-linear. It also
can be noticed that the $2$ extremal parameters, corresponding to the
extrapolation regime and represented in the plot by the two most outer
trajectories that enclose the other $6$, have smooth latent dynamics
analogously to the others even though the reconstruction error starts
degrading, as can be seen from Figure~\ref{fig: ml_ncl}.

\subsubsection{Hyper-reduction and teacher-student training}
The selection of the magic points is carried out with the greedy
Algorithm~\ref{algo:nodesRoc}. The FOM snapshots employed correspond to the training
parameters $1$ and $12$, but are sampled every $10$ time step instead of every
$4$, as for the training snapshots of the CAE.

We perform a convergence study increasing the number of magic points from $50$
to $100$ and $150$. The corresponding submesh sizes, i.e. the number of cells involved in the
discretization of the residuals, are reported in the Table~\ref{tab: decoder ts
    ncl}. The submesh size is bounded above with the total number of
the cells in the mesh, that is $3600$.

After the computation of the magic points, the FOM snapshots are restricted to
those cells and employed as training outputs of the compressed decoder, as described in
Section~\ref{subsec: compressed decoder training}. The actual dimension of the
outputs is twice the submesh size, since for each cell there are $2$ degrees of
freedom corresponding to the velocity components. The inputs are the
$4$-dimensional latent coordinates of the encoded FOM training snapshots, for a
total of $6012$ training input-output pairs. The architecture of the compressed
decoder is a feedforward neural network (FNN) with one hidden layer, whose number
of nodes is reported in Table~\ref{tab: decoder ts ncl}, under 'HL size'. The
compressed decoders architecture's specifics are also summarized in Table~\ref{tab: dec ncl}.

\begin{figure}[htp!]
    \centering
    \includegraphics[width=0.8\textwidth]{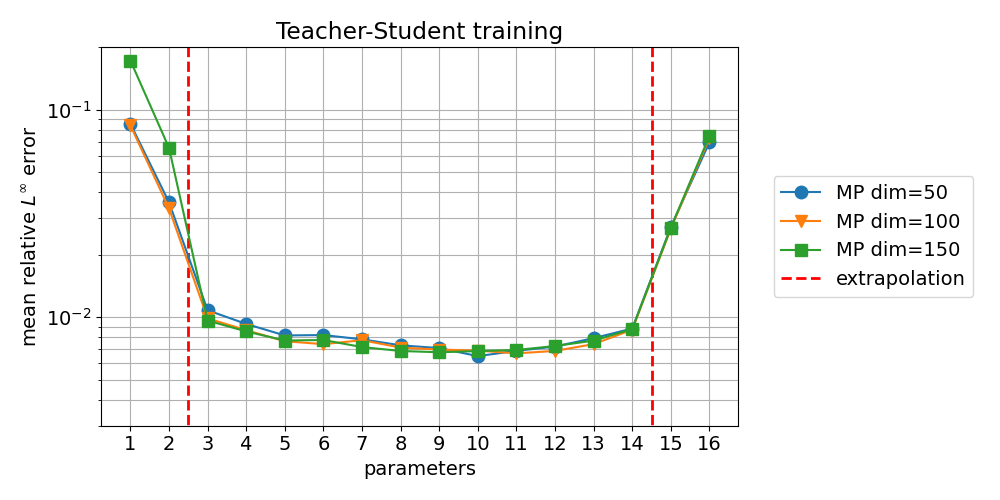}
    \caption{Prediction accuracy on the test snapshots restricted to the magic
        points. The accuracy is measured with the relative $L^{\infty}$-error
        averaged over the time trajectories associated to each one of the $N^{\mu}_{\text{test}}=16$ test samples.}
    \label{fig:decoder compress accuracy}
\end{figure}

Each compressed decoder is trained for $3000$ epochs, with an initial learning
rate of $1\mathrm{e}-4$, that halves if the loss from Equation~\ref{eq: decoder
    compress} does not decrease after $200$ epochs. The batch size is $20$. The duration of the training is
reported in Table~\ref{tab: decoder ts ncl} under 'TS total epochs', that stands
for Teacher-Student training total epochs, along with the average cost for an
epoch, under 'TS avg epoch'. The accuracy of the predictions on the test snapshots restricted to the magic
points is assessed in Figure~\ref{fig:decoder compress accuracy}.

\begin{table}[htp!]
    \centering
    \caption{In this table are reported for the NCL test case: the submesh size
        and the hidden layer (HL) number of nodes of the compressed decoder; the
        Teacher-Student training (TS) duration in secods (TS total epochs), the Teacher-Student
        training average epoch duration in seconds (TS avg epoch); the average time step for
        NM-LSPG-GNAT with compressed decoder (avg GNAT-TS) in milliseconds, and the
        average time step for NM-LSPG-GNAT with the hyper-reduced residuals but full
        CAE decoder (avg GNAT-no-TS) in milliseconds. In red, the results are obtained with $13$ maximum
        residual evaluations of the Levenberg-Marquardt algorithm, instead of the
        fixed $7$, see Remark~\ref{rmk: lm parameters}.}
    \footnotesize
    \begin{tabular}{ l | c | c | c | c | c | c}
        \hline
        \hline
        MP  & Submesh size & HL size & TS total epochs & TS avg epoch & avg GNAT-TS           & avg GNAT-no-TS \\
        \hline
        \hline
        50  & 139          & 300     & 1682 s          & 0.560 s      & 1.88 ms               & 4.496 ms       \\
        \hline
        100 & 246          & 350     & 2482 s          & 0.827 s      & \color{red}{4.599} ms & 4.499 ms       \\
        \hline
        150 & 335          & 400     & 4245 s          & 1.415 s      & 2.318 ms              & 4.208 ms       \\
        \hline
        \hline
    \end{tabular}
    \label{tab: decoder ts ncl}
\end{table}

\begin{figure}[htp!]
    \centering
    \includegraphics[width=1\textwidth]{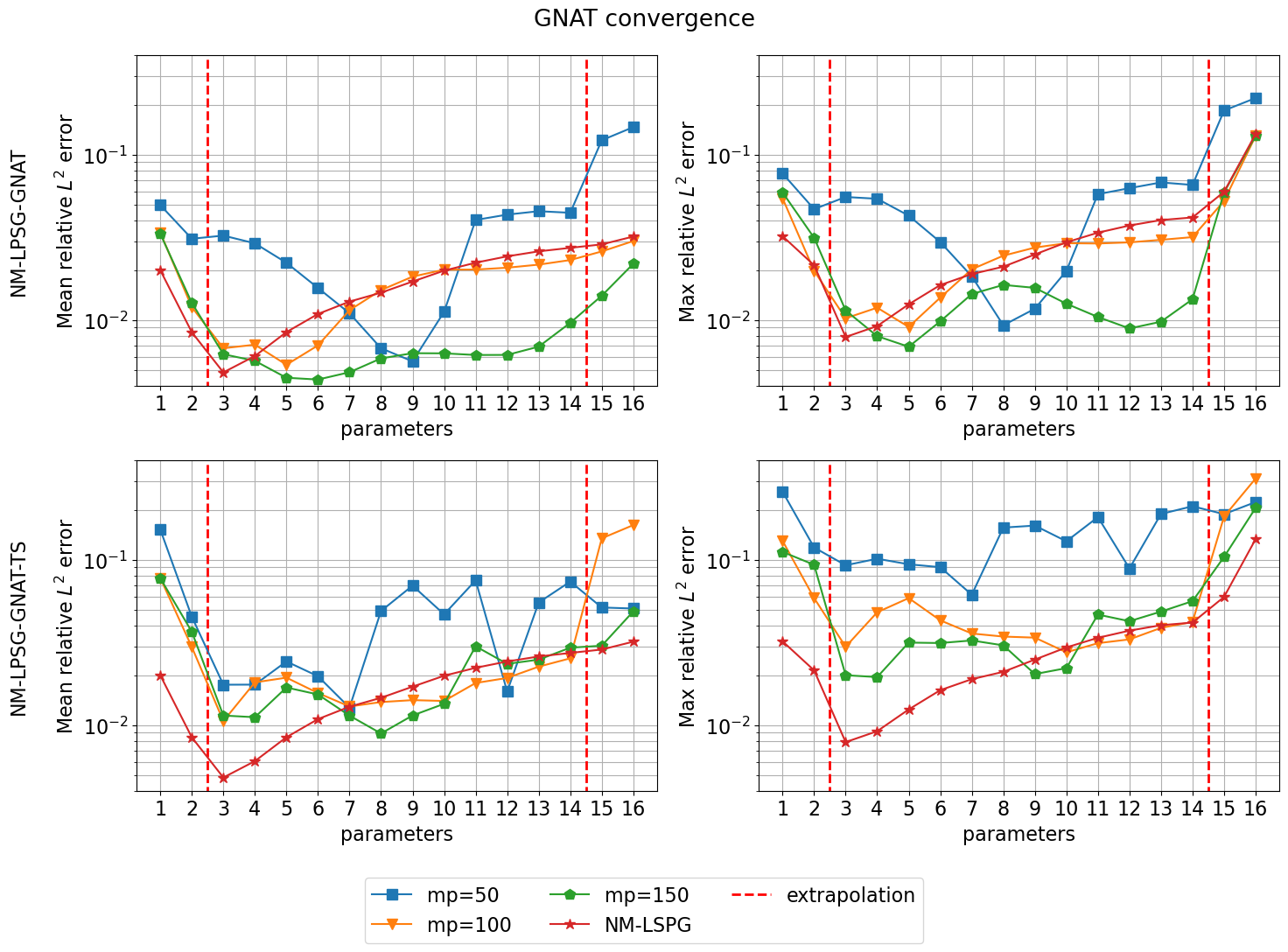}
    \caption{First row: NM-LSPG-GNAT mean and max relative $L^2$-error. Second row: NM-LSPG-GNAT-TS mean and max relative $L^2$-error. In red also the NM-LSPG accuracy is reported. The mean and max values are evaluated with respect to the time series of intermediate solutions associated to each one of the $N^{\mu}_{\text{test}}=16$ test parameters.}
    \label{fig: gnat convergence ncl}
\end{figure}

The convergence with respect to the number of magic points is shown in
Figure~\ref{fig: gnat convergence ncl}. Since, especially for the NM-LSPG-GNAT-TS
reduced-order model, the relative $L^2$-error is not
uniform along the time scale, we report both the mean and max relative
$L^2$-errors over the time series associated to each one of the $16$ test
parameters.

From Figure~\ref{fig: gnat convergence ncl} and Table~\ref{tab: decoder ts ncl}
it can be seen that NM-LPSG-GNAT is more accurate than NM-LSPG-GNAT-TS even
tough computationally more costly in the online stage. We underline that in the
offline stage NM-LSPG-GNAT-TS requries the training of the compressed decoder.
However, this could be performed at the same time of the CAE training, see the
discussion section~\ref{sec:discussion}. The NM-LSPG-GNAT reduced-order model
achieves better results also in the extrapolation error, sometimes even lower
than the NM-LSPG method: this remains true even when increasing the maximum residual
evaluations of the LM algorithm and it may be related to the nonlinearity of the decoder that
introduces difficult to interpret correlations of the latent dynamics with the
output solutions restricted to the magic points.

All the simulations of NM-LSPG, NM-LSPG-GNAT-TS and NM-LSPG-GNAT methods
converge with a maximum of $7$ residual evaluations of the LM algorithm, for all magic points
reported, that is $50$, $100$, and $150$ and for all the $16$ test parameters.
However, $3$ test parameters could not converge for the method NM-LSPG-GNAT-TS
with $100$ magic points, so we increased the maximum function evaluations to
$13$ for all the test points. The higher computational cost per time step is shown
in Table~\ref{tab: decoder ts ncl}. A part from those $3$ test points not
converging, the accuracy remains the same for the other $13$ test parameters, so
we have chosen to report the results in the case of $13$ residual evaluations
for all the $16$ test parameters in Figure~\ref{fig: gnat convergence ncl}.

\subsubsection{Comparison with data-driven predictions based on a LSTM}
To asses the quality of the reduced-order models devised, we compare the
accuracy in the training and extrapolation regimes, and the computational cost
of the offline and online stages with a purely data-driven ROM in which the
solutions manifold is approximated by the same CAE, but the dynamics is evolved
in time with a LSTM neural network. The architecture of the LSTM employed is
reported in Table~\ref{tab: lstm}. The results are summarized in Figure~\ref{fig:
    lstm comapre ncl}, the computational costs in Table~\ref{tab: lstm compare ncl}.

The LSTM is trained for $10000$ epochs with the ADAM stochastic optimization
algorithm and an initial learning rate of $0.001$, halved if after $500$ epochs
the loss does not decrease. The time series used for the training are the same
$N_{\text{train}}=6012$ training snapshots employed for the CAE. We remark that the LSTM cannot
approximate the dynamics for an arbitrary time step, but it is fixed, depending
on the training time step used, in this case $0.004$ seconds.

The offline stage's computational cost is determined by the heavy CAE training for both
the procedures, see the Discussion
section~\ref{sec:discussion} for possible remedies. The LSTM-NN achives a
speed-up close to $3$ with respect to the FOM, differently from the NM-LSPG-GNAT
and NM-LSPG-GNAT-TS methods. However, since the models are hyper-reduced,
increasing the degrees of freedom refining the mesh should increase the
computational cost of the FOM and NM-LSPG methods only, with due precautions.
The average cost of the evaluation of the dynamics for the LSTM model with a time step of $0.004$
seconds is associated to the label 'avg LSTM-NN full-dynamics'; thanks to
vectorization, the dynamics for all the $N^{\mu}_{\text{test}}=16$ parameters is evlauated with a
single forward of the LSTM, thus the low computational cost reported.

\begin{table}[htp!]
    \centering
    \caption{Offline stage: full-order model (FOM) computation of the
    $N_{\text{train}}=6012$ snapshots (FOM snapshots evaluation), average cost of a single epoch for the
    training of the CAE with a batch size of $20$ (CAE training single epoch avg), and total cost for $3000$
    epochs (CAE training); average epoch's cost for the LSTM training with a batch size of
    $100$, and total cost for $10000$ epochs. Online stage: for the FOM, NM-LSPG
    and LSTM-NN models it is reported the average of a single time step cost
    over all the $N_{\text{test}}=1616$ test parameters and time series, and the
    average cost of the full dynamics over the $N^{\mu}_{\text{test}}=16$ test
    parameters. The LSTM-NN full-dynamics is evaluated with a single forward.}
    \footnotesize
    \begin{tabular}{ l | c}
        \hline
        \hline
        Offline stage                     & Time        \\
        \hline
        \hline
        FOM snapshots evaluation          & 29.04 [s]   \\
        CAE training single epoch avg     & 10.1 [s]    \\
        CAE training                      & 20213.3 [s] \\
        \hline
        LSTM-NN training single epoch avg & 0.139 [s]   \\
        LSTM-NN training                  & 1365 [s]    \\
        \hline
        \hline
    \end{tabular}
    \hspace{1mm}
    \begin{tabular}{ l | c}
        \hline
        \hline
        Online stage              & Time        \\
        \hline
        \hline
        avg FOM time step         & 1.210 [ms]  \\
        avg FOM full dynamics     & 2.42 [s]    \\
        \hline
        avg NM-LSPG time step     & 22.126 [ms] \\
        avg NM-LSPG full-dynamics & 133.2 [s]   \\
        \hline
        avg LSTM-NN time step     & 0.432 [ms]  \\
        LSTM-NN full-dynamics     & 6.982 [ms]  \\
        \hline
        \hline
    \end{tabular}
    \label{tab: lstm compare ncl}
\end{table}

\begin{figure}[htp!]
    \centering
    \includegraphics[width=1\textwidth]{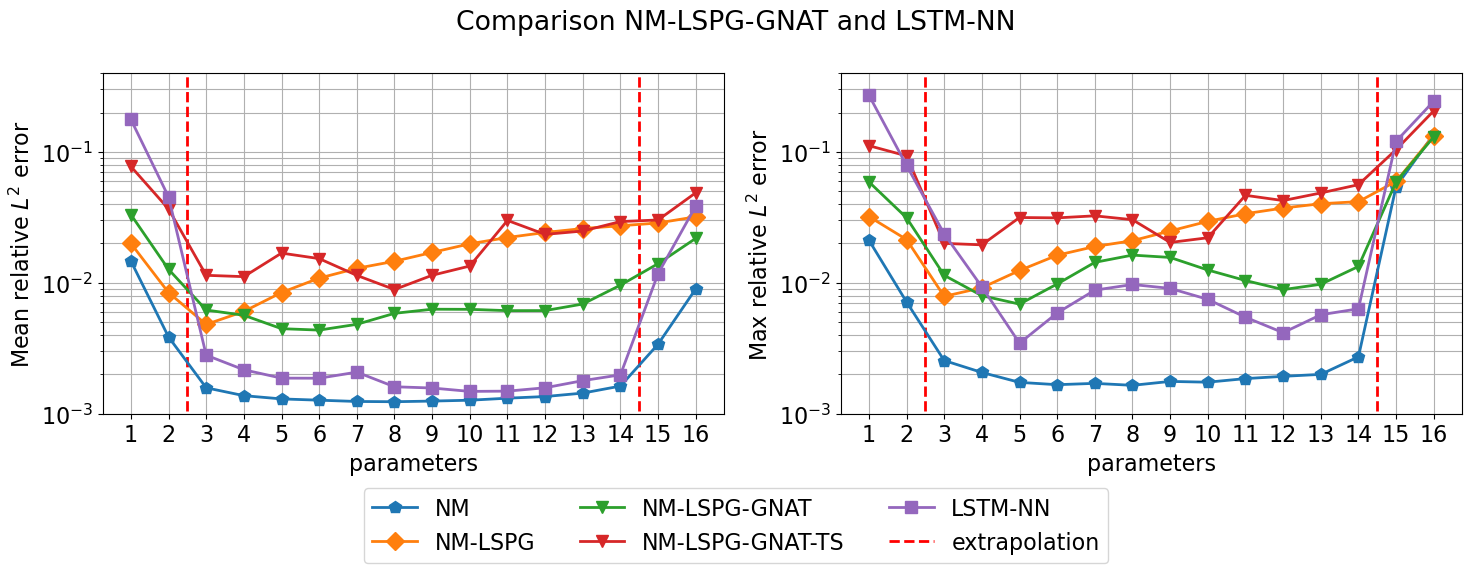}
    \caption{Comparison of the accuracies beetween all the ROMs presented on the test set of $N^{\mu}_{\text{test}}=16$ parameters, each associated to a time series of $N^{t}_{\text{test}}=101$ intermediate solutions. The number of magic points employed for NM-LSPG-GNAT and NM-LSPG-GNAT-TS is $150$. The mean and maximum relative $L^2$-errors are taken with respect to the time scale.}
    \label{fig: lstm comapre ncl}
\end{figure}

The CAE reconstruction error in blue in Figure~\ref{fig: lstm comapre ncl},
lower bounds all the other models' errors. This is the reason why having a good
accuracy of the CAE's solution manifold approximation is mandatory to build up
non-linear manifold methods. In this sense NM-LSPG-ROC-TS and NM-LSPG-GNAT-TS
with respect to NM-LSPG-GNAT with shallow autoencoders~\cite{kim2020fast} offer the possibility
to choose an arbitrary architecture for the autoencoder, thus allowing a more
accurate solution manifold approximation.

While in the training range from test parameter $3$ to $14$, the accuracy of the
LSTM-NN is significantly better than NM-LSPG-GNAT and NM-LSPG-GNAT-TS models',
in the extrapolation regime we observe that the predictions of the fully
data-driven model degrades. The extrapolation error of the LSTM-NN model depends
on the architecture chosen, regularization applied, training procedure, and
hyper-parameters tuning. What
can be assessed from the results is that, outside the training range, the
LSTM-NN's accuracy is dependent on all these factors, with sometimes a difficult
interpretation of the results, while NM-LSPG-GNAT relies
only on the number of magic points employed and the dynamics is evolved in time
minimizing a physical residual directly related to the NCL model's equations.

\subsection{Shallow Water Equations (SWE)}
\label{subsec: SWE}
The second test case we present is a 2d non-linear, time-dependent, parametric model based on the shallow water equations
(SWE). Also in this case, the non-temporal parameter affects the initial
conditions, $\mu\in[0.1, 0.3]$:
\begin{equation}
    \label{eq:shallow 2d}
    \begin{cases}
        \partial_{t}(h\mathbf{u}) + \nabla\cdot(h\mathbf{u}\otimes\mathbf{u}) + \frac{1}{2}|\mathbf{g}|h\nabla h= 0                                                                                                        & (\mathbf{x}, t)\in[0, 1]^{2}\times [0, 0.2],        \\
        \partial_{t} h + \nabla\cdot(\mathbf{u}h) = 0                                                                                                                                                                      & (\mathbf{x}, t)\in[0, 1]^{2}\times [0, 0.2],        \\
        \mathbf{u}(\mathbf{x}, 0) = 0                                                                                                                                                                                      & \mathbf{x}\in[0, 1]^{2},                            \\
        h(\mathbf{x}, 0) = \mu\left(\frac{1}{e^{1}} \cdot e^{-\frac{1}{0.04 - \lVert\mathbf{x}-\mathcal{O}\rVert_{2}^{2}}}\chi_{\lVert\mathbf{x}\rVert_{2}^{2}<0.2} + \chi_{\lVert\mathbf{x}\rVert_{2}^{2}\geq 0.2}\right) & \mathbf{x}\in[0, 1]^{2},                            \\
        \mathbf{u}(\mathbf{x}, t)\cdot\mathbf{n} = 0                                                                                                                                                                       & (\mathbf{x}, t)\in\partial[0, 1]^{2}\times[0, 0.2], \\
        \nabla h(\mathbf{x}, t)\cdot\mathbf{n} = 0                                                                                                                                                                         & (\mathbf{x}, t)\in\partial[0, 1]^{2}\times[0, 0.2],
    \end{cases}
\end{equation}
where $h$ is the water depth, $\mathbf{u}$ is the velocity vector,
$\mathbf{g}$ is the gravitational acceleration, and $\mathcal{O}$ is the point $(0.5,
    0.5)\in\Omega$. We consider a constant batimetry $h_0=0$, so that the free
surface height $h_{\text{total}}=h+h_0$ is equal to the water depth $h$.

The time step that will be employed for the evolution of the dynamics of the FOM
is $1\mathrm{e}-4$ seconds. The training and test snapshots are sampled every $4$ time
steps. The training non-temporal parameters are $N^{\mu}_{\text{train}}=10$ in number and they are
sampled equispacedly in the training interval $\mu\in[0.1, 0.3]$, for a total of
$N_{\text{train}}=N^{\mu}_{\text{train}}\cdot N^{t}_{\text{train}} = 10\cdot
501=5010$ training snapshots.

Due to an inaccurate reconstruction
erorr of the CAE for the first time instants, the predictions of the dynamics of the reduced model are
evaluated from the time instant $t_0 = 0.01$ seconds. The test
non-temporal parameters are $N^{\mu}_{\text{test}}=8$ in number and sampled equispacedly in the test
interval $\mu\in[0.05, 0.35]$, for a total of
$N_{\text{test}}=N^{\mu}_{\text{test}}\cdot N^{t}_{\text{test}} = 8\cdot
475=3800$ test snapshots,
since the first $26$ are cut from the time series, $N^{t}_{\text{train}}=501$ snapshots long.
Again the first $2$ and the last $2$ parameters
correspond to the extrapolation regime. The initial latent
variables are obtained projecting with the encoder into the latent space the
test snapshots corresponding to the time instants $t_0=0.01$ instead of $t=0$. The training and test time series are
labelled from $1$ to $10$ and from $1$ to $8$ with an increasing order.

\label{subs:shallow water}
\begin{figure}[h]
    \includegraphics[width=0.24\textwidth]{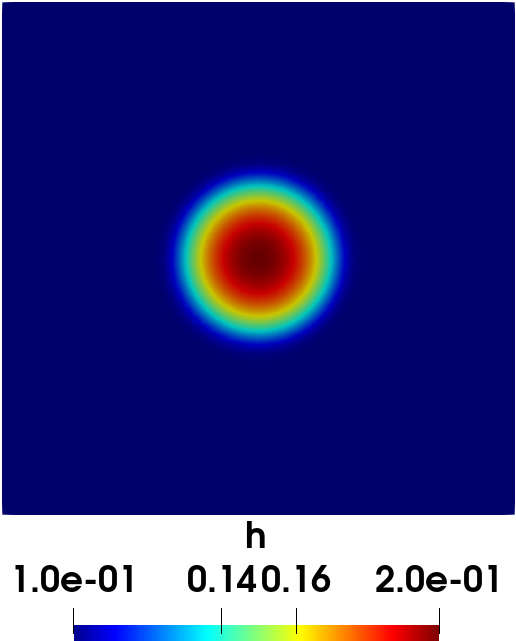}
    \includegraphics[width=0.24\textwidth]{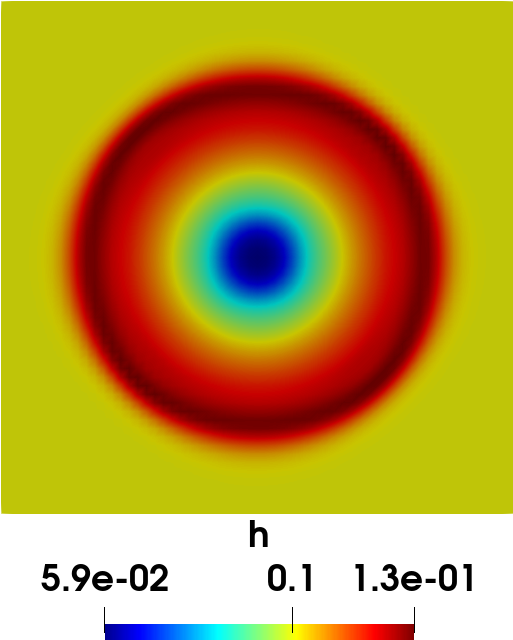}
    \includegraphics[width=0.24\textwidth]{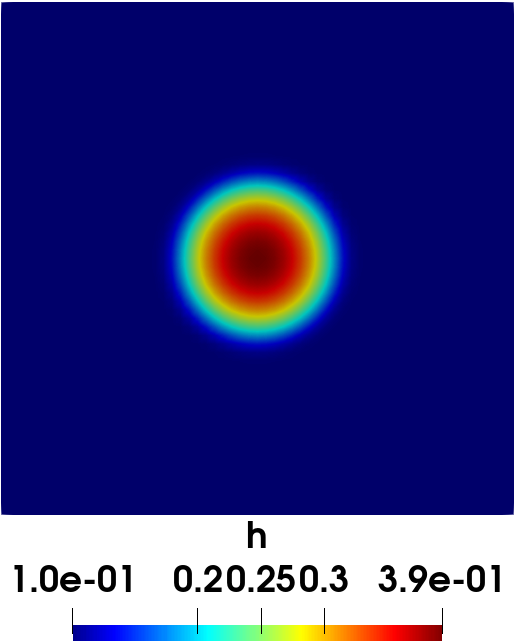}
    \includegraphics[width=0.24\textwidth]{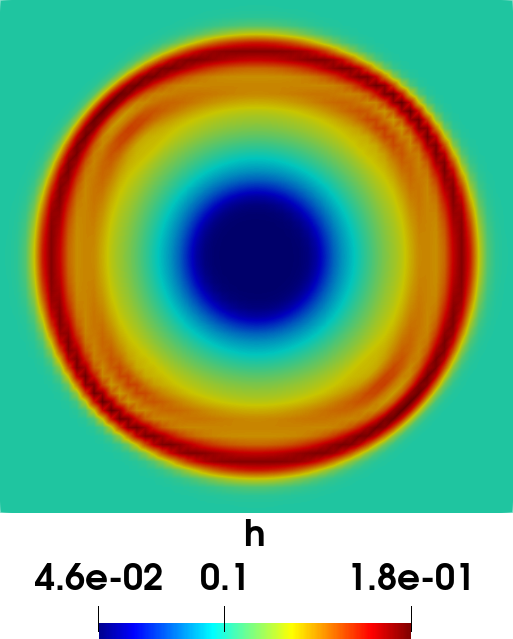}
    \caption{From left to right: FOM solution of equation~\ref{eq:burgers 2d} at
        $(t, \mu)\in\{(0.01, 0.1), (0.2, 0.1), (0.01, 0.3),(0.2, 0.3)\}$.}
    \label{fig:qualitative snap sw h}
\end{figure}

\begin{figure}[h]
    \includegraphics[width=0.24\textwidth]{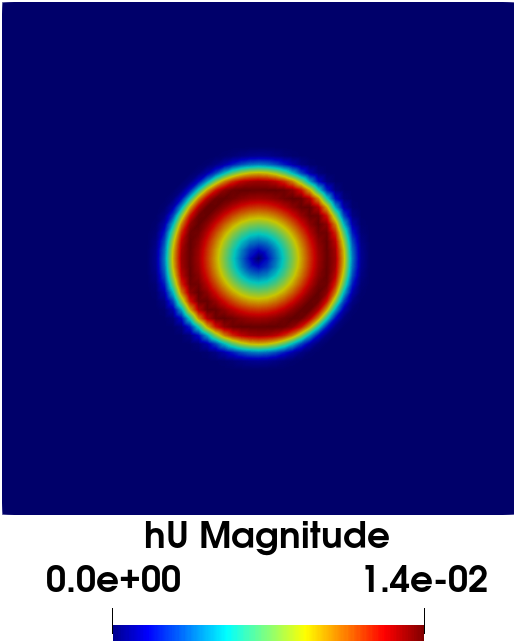}
    \includegraphics[width=0.24\textwidth]{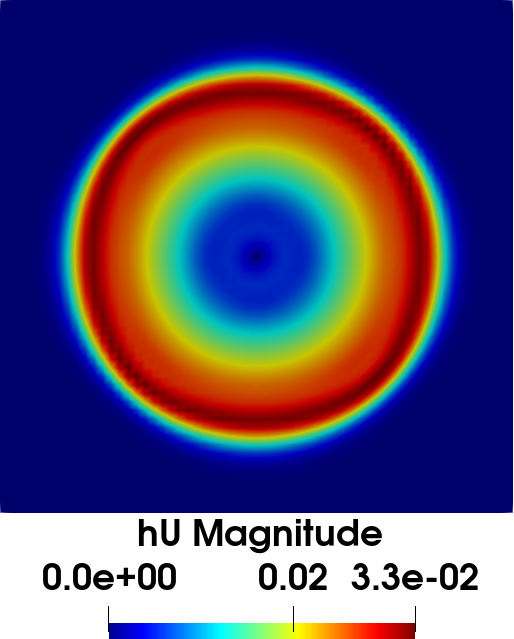}
    \includegraphics[width=0.24\textwidth]{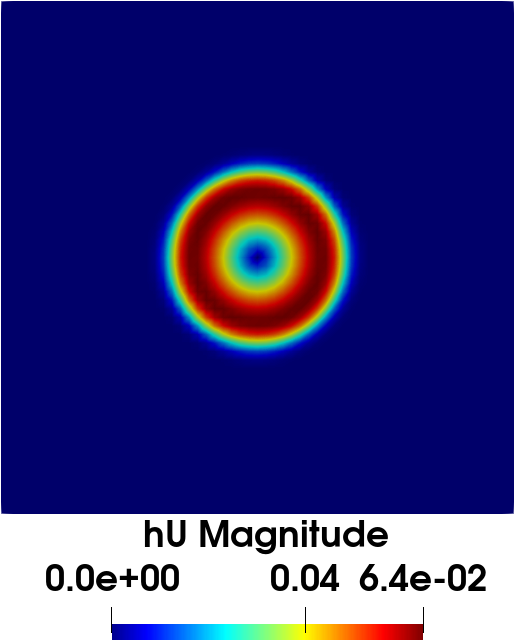}
    \includegraphics[width=0.24\textwidth]{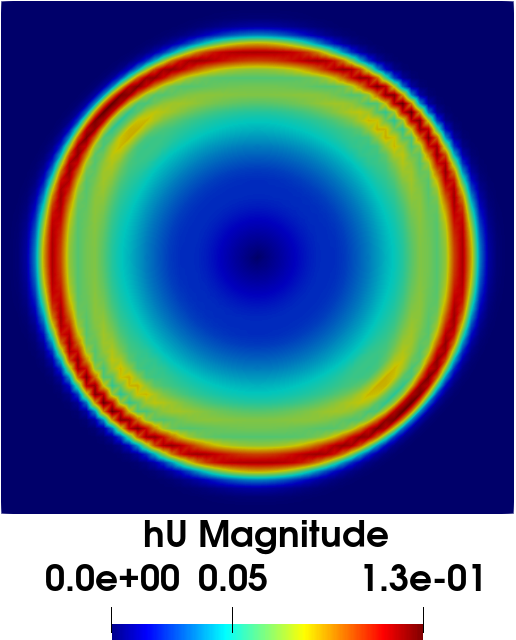}
    \caption{From left to right: FOM solution of equation~\ref{eq:burgers 2d} at
        $(t, \mu)\in\{(.010, 0.1), (0.2, 0.1), (0.01, 0.3),(0.2, 0.3)\}$.}
    \label{fig:qualitative snap sw hU}
\end{figure}

In this test case the ROC hyper-reduction is more accurate with respect to the
GNAT one, so the results are reported w.r.t. this hyper-reduction method.

\subsubsection{Full-order model}
One detail that we didn't stress in the previous test case is that the FOM and
the NM-LSPG ROM can discretize the residuals of the SWE differently: only the
consistency of the discretization is required, characterizing the NM-LPSG family
as equations-based rather than fully intrusive.

The FOM solutions are computed with the OpenFoam solver
\textit{shallowWaterFoam}~\cite{OF}, while the ROMs discretize the residual with
a much simpler numerical scheme.

The FOM numerical scheme is the PIMPLE algorithm, a combination of
PISO~\cite{patankar1983calculation} (Pressure Implicit with Splitting of
Operator) and SIMPLE~\cite{issa1991solution} (Semi-Implicit Method for
Pressure-Linked Equations). For the shallow water equations the free surface
height $h$
plays the role of the pressure in the Navier-Stokes equations, regarding the
PIMPLE algorithm implementation. The number of outer PISO corrections is $3$.

The time discretization is performed with the semi-implicit Euler method. The
non-linear advection terms are discretized with the Linear-Upwind Stabilised
Transport (LUST) scheme that requires a stencil with $2$ layers of adjacent
cells for the hyper-reduction. The gradients are linearly interpolated through
Gauss formula. The solutions for $hU$ are obtained with Gauss-Seidel iterative
method, and for $h$ with the conjugate gradient method preconditioned by the
Diagonal-based Incomplete Cholesky (DIC) preconditioner. For both of them the absolute tolerance on
the
residual is $1\mathrm{e}-6$ and the relative tolerance of the residual w.r.t.
the initial condition is $0.1$.

The residual of ROMs is instead discretized as follows. If we represent
with $M_{hU},\ M_h$ the mass matrices, with $G(h)$ the discrete gradient vector of $h$,
and with $C_{hU}((hU)^{t-1}),\ C_{h}(U^{t-1})$ the
advection matrices, then, at every time
instant $t$, the discrete equations
\begin{align}
    \frac{M_{hU}}{\Delta t}(hU)^t+C_{hU}((hU)^{t-1})U^t + gh G(h)= \frac{M_{hU}}{\Delta t}(hU)^{t-1}, \\
    \frac{M_{h}}{\Delta t}h^t+C_h(U^{t-1})h^t  = \frac{M_{h}}{\Delta t}h^{t-1},
\end{align}
are solved for the state $((hU)^t, h^t)$ with a semi-implicit Euler method. The
same numerical schemes and linear systems iterative solvers of the FOM are
employed. In principle, they could be changed.

Since now the stencil of a single cell needs two layers of adjacent cells for
the discretizations, for each internal magic point, $12$ additional cells need
to be considered for the hyper-reduction.

\subsubsection{Manifold Learning}
\label{subsec: ml sew}
\begin{figure}[htp!]
    \centering
    \includegraphics[width=0.9\textwidth]{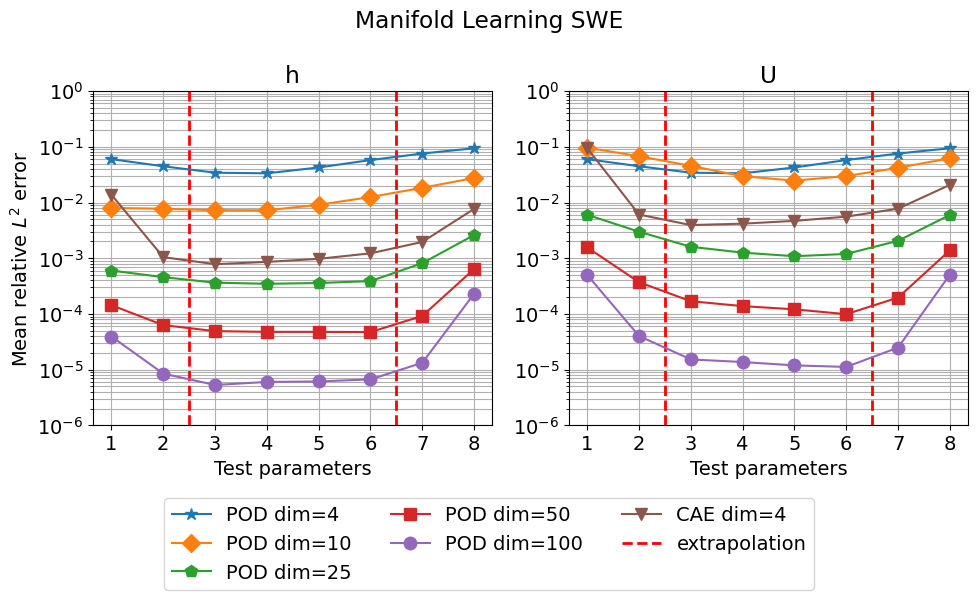}
    \caption{Comparison between the CAE's and POD's projection errors of the discrete solution manifold represented by the $N^{\text{train}}=N^{\mu}_{\text{train}}\cdot N^{t}_{\text{train}}=10\cdot 501=5010$ training snapshots. The error is evaluated on the $8$ test parameters, each associated with a time series of $N^{t}_{\text{test}}=475$, for a total of $N_{\text{test}}=3800$ test snapshots. The mean is performed over the time scale.}
    \label{fig: nm swe}
\end{figure}

The CAE architecure for the SWE model is reported in Table~\ref{tab: cae swe}. This
time one encoder and two decoders, one for the velocity $U$ and one for the
height $h$ are trained. Moreover, to increase the generalization capabilities we
converted $2$ layers of the decoder for $U$ in recurrent convolutional layers
as shown in the Appendix. Even with this modification the initial time steps,
from $0$ to $0.01$ are associated to a high reconstruction error: the relative $L^2$-error
is around $0.1$ for every test parameter at the initial time instants, slowly
decreasing towards the accuracy shown in Figure~\ref{fig: nm swe} after
$t=0.01$, chosen as initial instant from here onward.

The CAE is trained for $500$ epochs with a batch size of $20$ and an initial
learning rate of $1\mathrm{e}-4$, that halves if the loss from
Equation~\ref{eq:loss cae} does not decrease after $50$ epochs. In this case the
state is $(U, h)$ so the encoder has three channels, two for the velocity
components, $U_1$, $U_2$, and one for the height, $h$. The high computational cost is shown in
Table~\ref{tab: swe offline online costs}. We have to observe that nor the
architecture is parsimonious for a good approximation of the discrete solution
manifold in the time interval $[0.01, 0.2]$, neither the number of training snapshots
$5010$ is optimized to reach the highest efficiency with the lowest
computational cost. Our focus is obtaining a satisfactory reconstruction error
in order to build up our ROMs.

The solution manifold parametrized by the decoder achieves the reconstruction
error of a linear manifold spanned by around $20$ POD modes. In fact, the decay
of the reconstruction error associated to the POD approximations is faster than
the previous test case in the time interval $[0.01, 0.2]$.

\begin{figure}[htp!]
    \centering
    \includegraphics[width=1\textwidth]{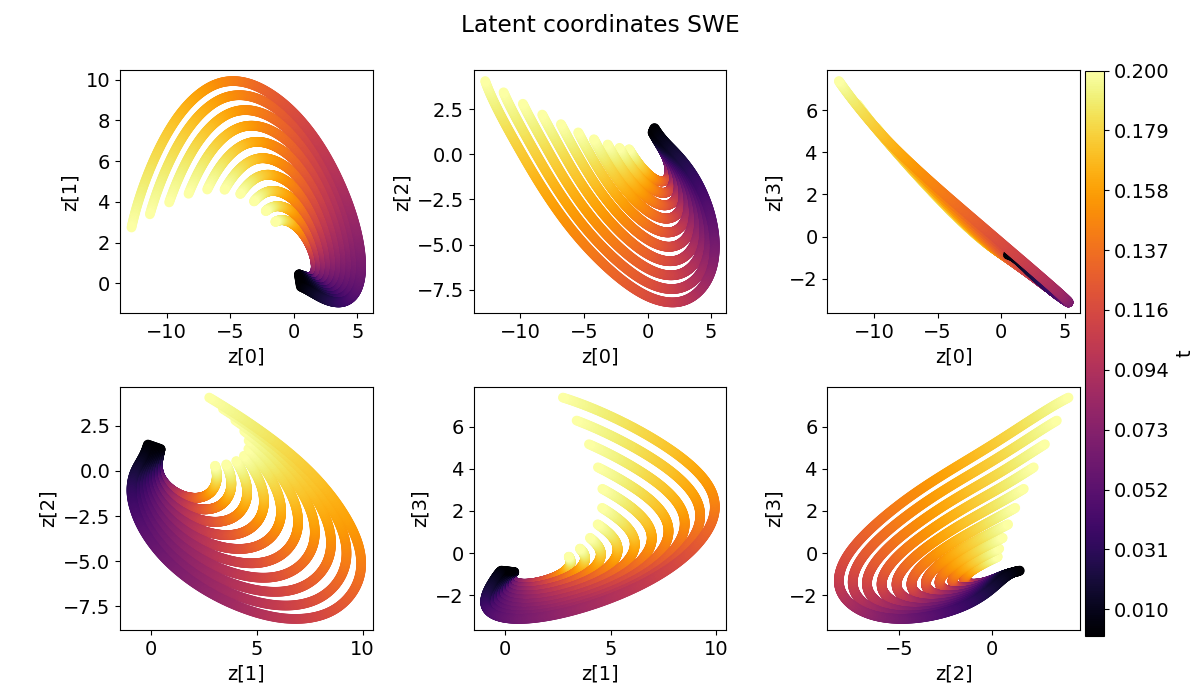}
    \caption{Correlations among the $4$ latent coordinates of the train set, $10$
        in total. They are obtained projecting the train snapshots into the latent
        space of the CAE with the encoder. The coloring corresponds to the time
        instants from $0.$ to $0.2$, with a time step of $4\mathrm{e}-4$ seconds.
        The initial time steps employed in the ROMs is $0.01$ seconds.}
    \label{fig: latent_test_swe}
\end{figure}

It can be seen from the representation of the latent dynamics associated to the
train parameters in Figure~\ref{fig: latent_test_swe}, that the initial
solutions overlap. This and
the low accuracy could
be explained by the fact that the FOM dynamics has different scales, especially for
the velocity $U$ that from the initial constant zero solution reaches a
magnitude of $10^{-1}$ meter per seconds. Further observations and possible
solutions are presented in the Discussion section~\ref{sec:discussion}.

\subsubsection{Hyper-reduction and teacher-student training}
\label{subsec: hr swe}
Mimicking the structure of the CAE, the compressed decoder is split in two, one
for the velocity $U$ and one for the free surface height $h$. The architecture
for both the decoders is a feed-forward NN with a single hidden layer; they are
reported in the Table~\ref{tab: dec swe}. The compressed decoders are trained for $1500$ epochs, with a batchsize of $20$,
and an initial learning rate of $1\mathrm{e}-4$ that halves after $100$ epochs if
the loss does not decreases, with a minimum value of $1\mathrm{e}-6$. As for the previous test case, the number of magic points, hidden layer sizes,
training times, average training epoch computational cost are reported in
Table~\ref{tab: decoder costs swe}. This time for each magic point correspond
$3$ degrees of freedom, so the actual output dimension of the compressed
decoders is three times the submesh sizes.

\begin{figure}[htp!]
    \centering
    \includegraphics[width=0.8\textwidth]{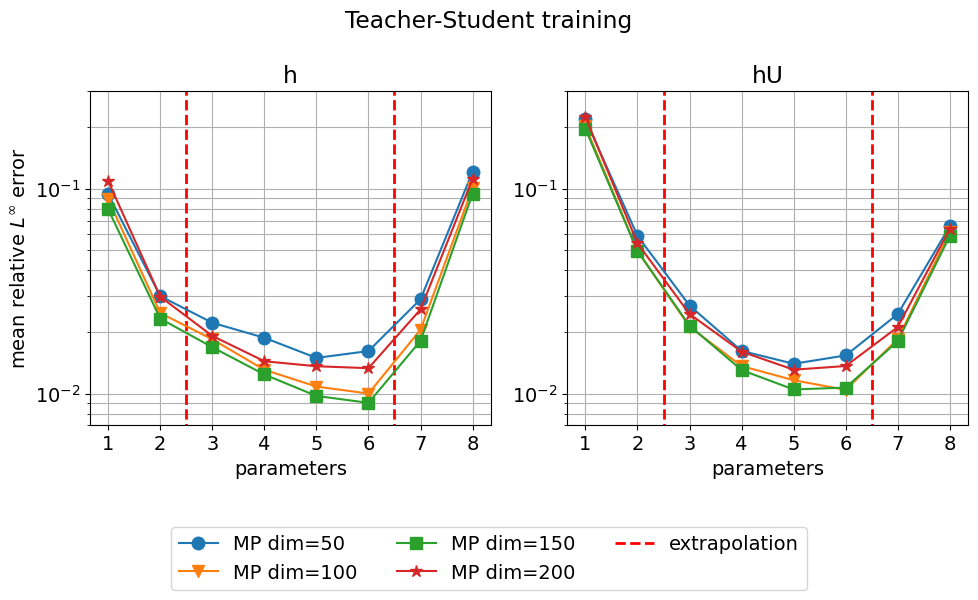}
    \caption{Prediction accuracy on the test snapshots restricted to the magic
        points. The accuracy is measured with the relative $L^{\infty}$-error
        averaged over the time trajectories associated to each one of the $N^{\mu}_{\text{test}}=8$ test
        samples.}
    \label{fig:dec swe accuracy}
\end{figure}

The relative $L^{\infty}$-error of the compressed decoder is shown in
Figure~\ref{fig:dec swe accuracy}. This time the extrapolation error is sensibly
higher as already seen in the reconstruction error of the CAE.

\begin{table}[htp!]
    \centering
    \caption{In this table are reported for the SWE test case: the submesh size
        and the hidden layer (HL) number of nodes of the compressed decoder; the
        Teacher-Student training (TS) duration in secods, the Teacher-Student
        training average epoch duration in seconds; the average time step for
        NM-LSPG-GNAT with compressed decoder (GNAT-TS) in milliseconds, and the
        average time step for NM-LSPG-GNAT with the hyper-reduced residuals but full
        CAE decoder (GNAT-no-TS) in milliseconds. In red the average computational cost with
        $13$ maximum residual evaluations due to the instability in the evolution of
        the dynamics of the test parameter $4$.}
    \footnotesize
    \begin{tabular}{ l | c | c | c | c | c | c}
        \hline
        \hline
        MP  & Submesh size & HL size & TS total epochs & TS avg epoch & avg GNAT-TS & avg GNAT-no-TS                    \\
        \hline
        \hline
        25  & 238          & 600     & 582 [s]         & 0.388432 [s] & 3.227 [ms]  & 14.782 [ms]                       \\
        \hline
        50  & 345          & 600     & 1234 [s]        & 0.823084 [s] & 4.062 [ms]  & 13.742 [ms]                       \\
        \hline
        100 & 590          & 900     & 6428 [s]        & 4.285968 [s] & 4.387 [ms]  & 15.238/{\color{red}{22.688}} [ms] \\
        \hline
        150 & 654          & 900     & 6746 [s]        & 4.497579 [s] & 4.307 [ms]  & 14.307 [ms]                       \\
        \hline
        \hline
    \end{tabular}
    \label{tab: decoder costs swe}
\end{table}

\begin{figure}[htp!]
    \centering
    \includegraphics[width=1\textwidth]{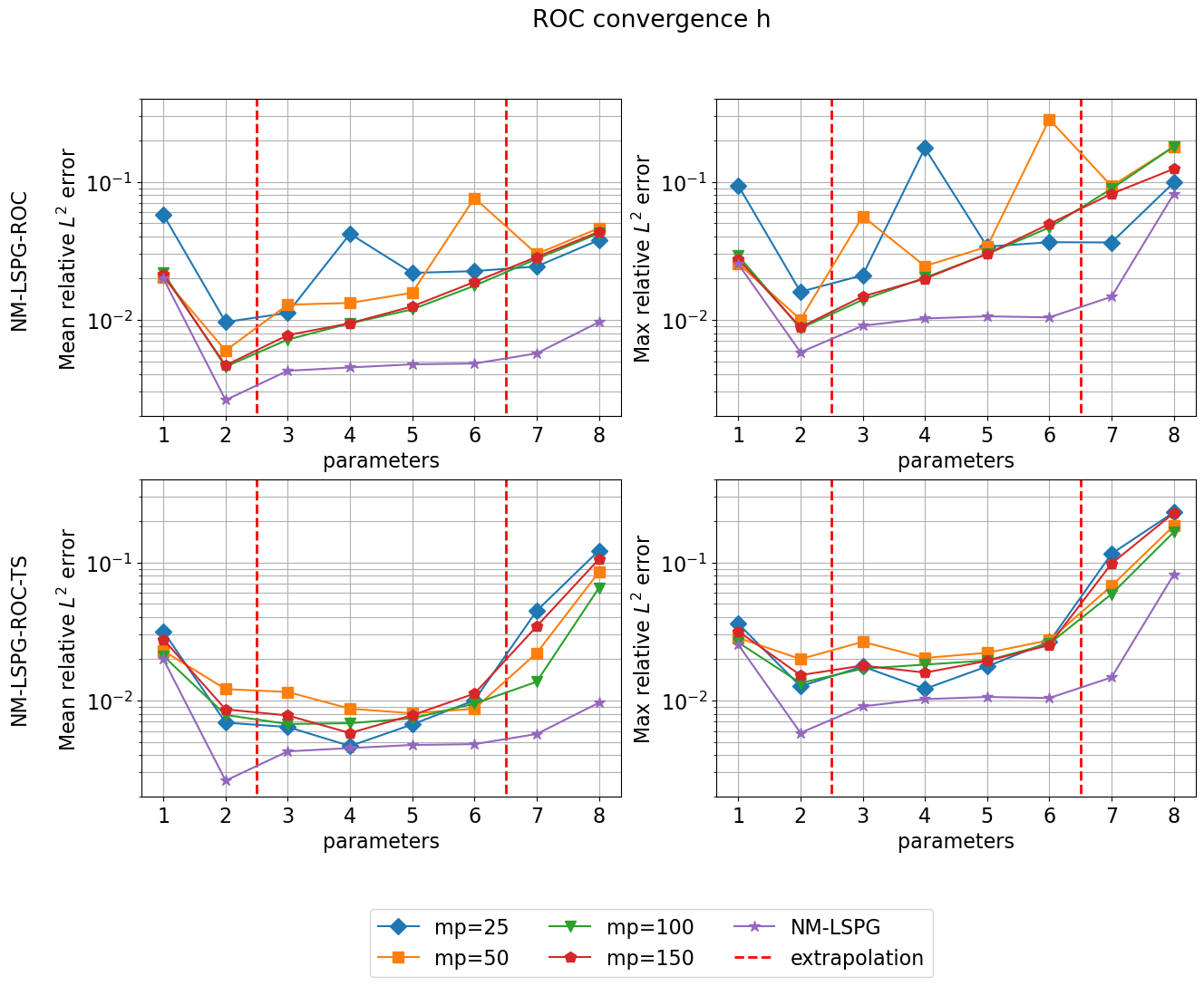}
    \caption{First row: NM-LSPG-ROC mean and max relative $L^2$-error. Second row: NM-LSPG-ROC-TS mean and max relative $L^2$-error. In violet also the NM-LSPG accuracy is reported. The mean and max values are evaluated with respect to the time series of intermediate solutions associated to each one of the $8$ test parameters.}
    \label{roc convergence swe h}
\end{figure}

As in the previous case, both the NM-LSPG-ROC and NM-LSPG-ROC-TS ROMs
are considered. This time it's the NM-LSPG-ROC's dynamics to be less stable with $7$
maximum residual evaluations of the LM algorithm. In this respect, for parameter
test $4$ and $\text{mp}=100$ the maximum number of residual evaluations is
increased to $13$; in the error plots in Figures~\ref{roc convergence swe h} and
\ref{roc
    convergence swe hU}, only the value for parameter $4$, $\text{mp}=100$ is
substituted. It is relevant to notice that the extrapolation error is lower for
higher numbers of magic points and for the NM-LSPG-ROC method.

\begin{figure}[htp!]
    \centering
    \includegraphics[width=1\textwidth]{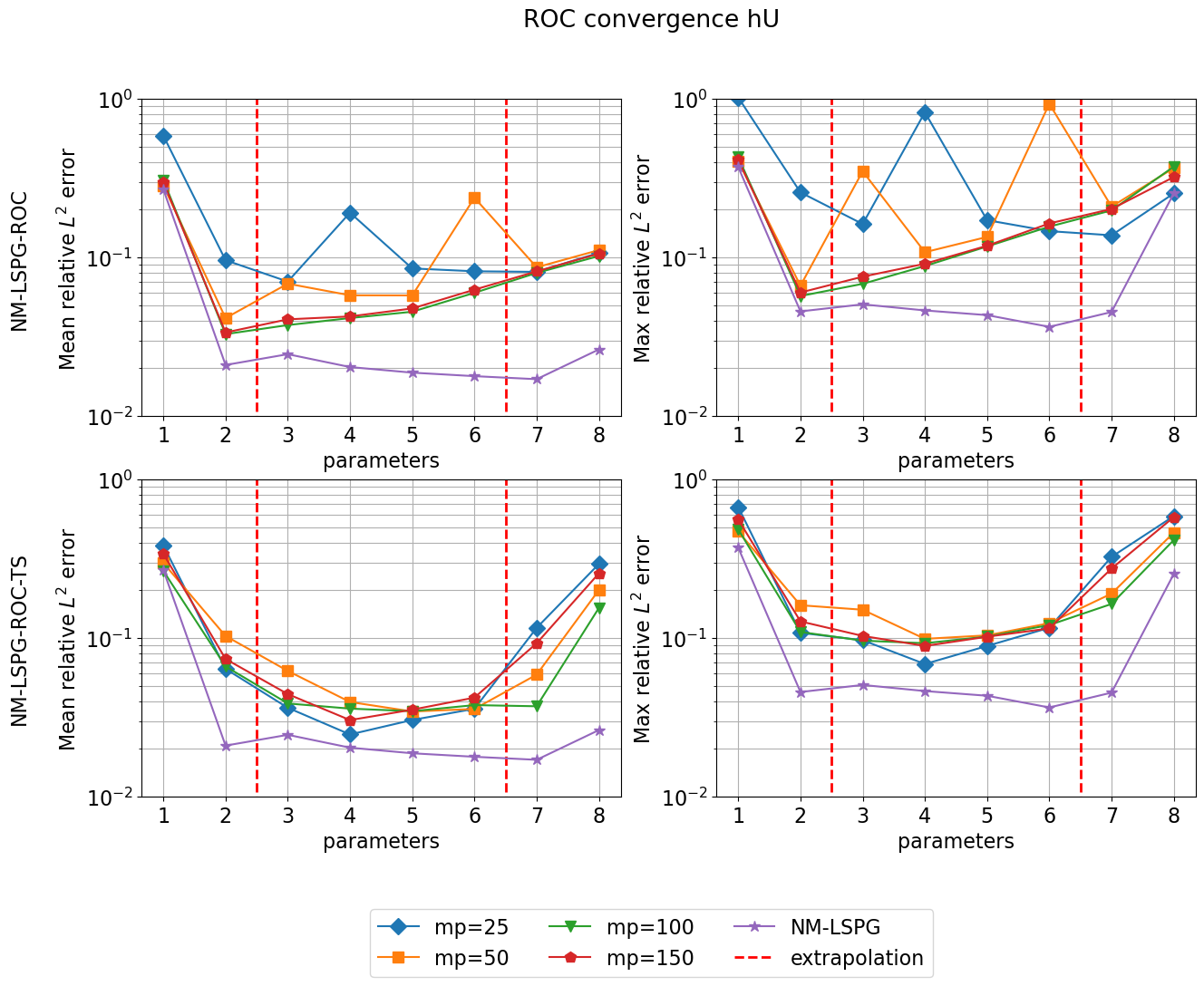}
    \caption{First row: NM-LSPG-ROC mean and max relative $L^2$-error. Second row: NM-LSPG-ROC-TS mean and max relative $L^2$-error. In violet also the NM-LSPG accuracy is reported. The mean and max values are evaluated with respect to the time series of intermediate solutions associated to each one of the $8$ test parameters.}
    \label{roc convergence swe hU}
\end{figure}

\subsubsection{Comparison with data-driven predictions based on LSTM}
The same LSTM architecture of the NCL test case is trained with the same
training procedure, only that now there are $5010$ training parameters-latent
coordinates pairs. For the architecture's specifics see Table~\ref{tab: lstm}. The computational costs introduced in the previous test case
are reported also for the SWE model in Table~\ref{tab: swe offline online
    costs}.

\begin{figure}[htp!]
    \centering
    \includegraphics[width=1\textwidth]{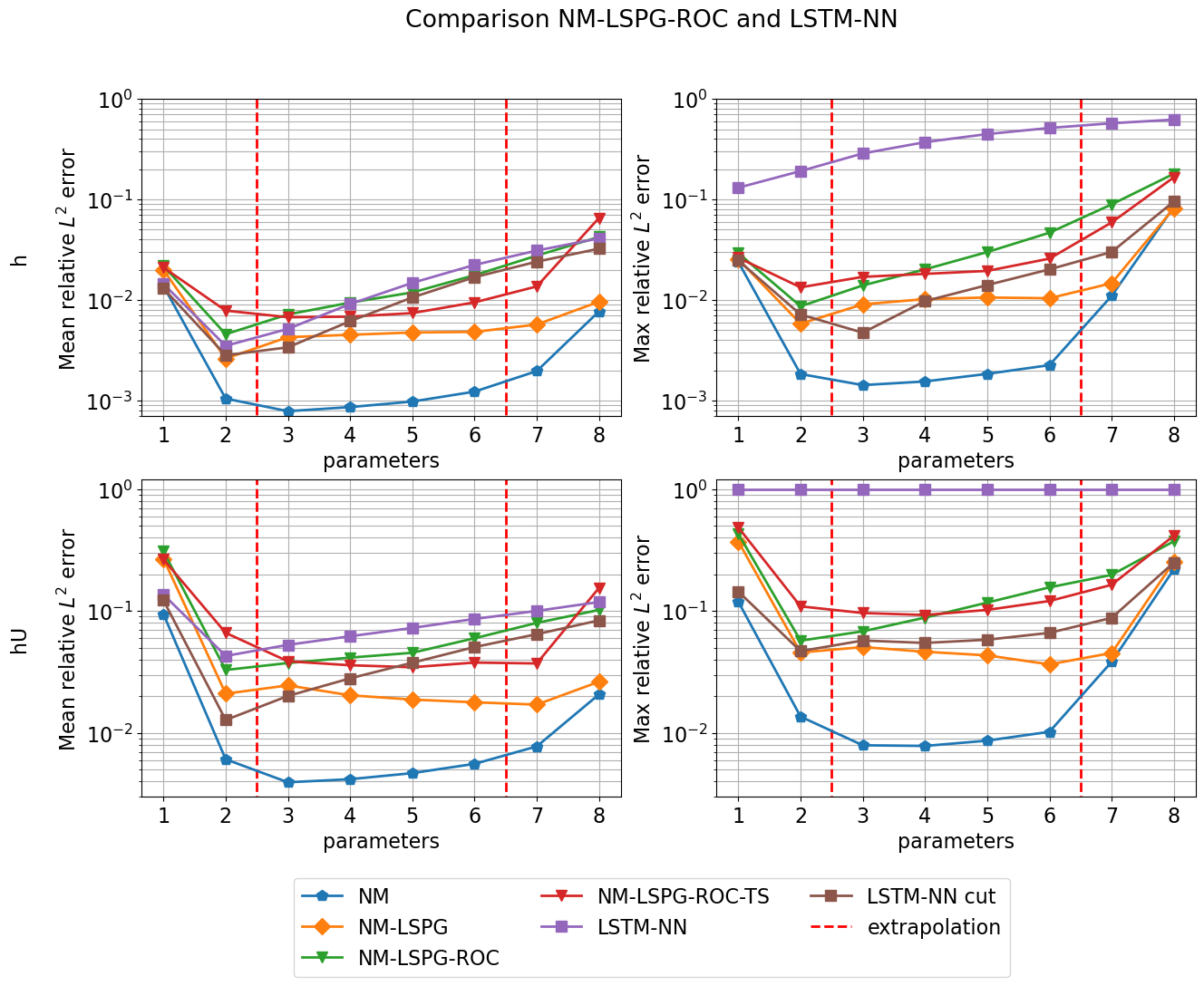}
    \caption{Comparison of the accuracies between all the ROMs presented on the
        test set of $8$ parameters, each associated to a time series of $475$
        intermediate solutions. The number of magic points employed for NM-LSPG-ROC
        and NM-LSPG-ROC-TS is $100$. The mean and maximum relative $L^2$-errors are
        taken with respect to the time scale. Since the LSTM is inaccurate between
        the time instants $0.01$ and $0.017$ seconds we report also the mean over
        the time interval $[0.017,\ 0.2]$ witch label 'LSTM-NN cut', to
        establish a fair comparison.}
    \label{lstm compare swe}
\end{figure}

With the architecture and training procedure employed, the LSTM could not
achieve a good accuracy at the initial time instants after $t_0=0.01$ seconds:
for this reason in the plot of the errorsini Figure~\ref{lstm compare swe} is reported
both the mean over the whole test time series of $475$ elements and over the
time series after the $40$-th element of the $475$, that corresponds to the the
time instant $0.017$ seconds. This issue should be ascribed at what we
discussed about in subsection~\ref{subsec: ml sew}, about the latent dynamics overlappings. Further remarks are provided
in the Discussion section~\ref{sec:discussion}.

A part from this, the LSTM predictions are for almost every test
parameter above only the NM-LSPG and CAE's reconstruction errors. Even in the
extrapolation regimes, the predictions are more accurate than the NM-LSPG-ROC
and NM-LSPG-ROC-TS reduced-order models. Regarding the computational costs, this time not
only the LSTM model but also the NM-LSPG-ROC-TS ROMs achieve a little speed-up
w.r.t. the FOM. The choice of doubling the decoders has repercussions in the
online costs, but it was made only in an effort to reach a good reconstruction
error of the CAE; maybe more light architectures could be employed for the
restricted time interval $[0.01, 0.2]$.

\begin{table}[htp!]
    \centering
    \caption{Offline stage: full-order model (FOM) computation of the $5010$ snapshots, average cost of a single epoch for the training of the CAE with a batch size of $20$, and total cost for $500$ epochs; average epoch's cost for the LSTM training with a batch size of $100$, and total cost for $10000$ epochs. Online stage: for the FOM, NM-LSPG and LSTM-NN models it is reported the average of a single time step cost over all the $8$ test parameters and time series, and the average cost of the full dynamics over the $8$ test parameters. The LSTM-NN full-dynamics is evaluated with a single forward.}
    \footnotesize
    \begin{tabular}{ l | c}
        \hline
        \hline
        Offline stage                     & Time           \\
        \hline
        \hline
        FOM snapshots evaluation          & 137.079881 [s] \\
        CAE training single epoch avg     & 94.700484 [s]  \\
        CAE training                      & 51431 [s]      \\
        \hline
        LSTM-NN training single epoch avg & 0.107722 [s]   \\
        LSTM-NN training                  & 1076.150 [s]   \\
        \hline
        \hline
    \end{tabular}
    \hspace{1mm}
    \begin{tabular}{ l | c}
        \hline
        \hline
        Online stage              & Time           \\
        \hline
        \hline
        avg FOM time step         & 7.251 [ms]     \\
        avg FOM full dynamics     & 116.024578 [s] \\
        \hline
        avg NM-LSPG time step     & 22.126 [ms]    \\
        avg NM-LSPG full-dynamics & 336.171916 [s] \\
        \hline
        avg LSTM-NN time step     & 1.802 [ms]     \\
        avg LSTM-NN full-dynamics & 6.834247 [s]   \\
        \hline
        \hline
    \end{tabular}
    \label{tab: swe offline online costs}
\end{table}

\section{Discussion}
\label{sec:discussion}
We comment the numerical results obtained:
\begin{itemize}
    \item Computational cost of the CAEs training. It is evident from
    Tables~\ref{tab: lstm compare ncl} and \ref{tab: dec ncl} that the offline
    stage's computational cost is dominated by the CAEs trainings. We have to
    remark that the architectures were not optimized to be the most persimonious
    ones in order to achieve the desired reconstruction error. Moreover,
    libtorch training took almost twice more time than the same architecture's
    training in PyTorch, due to implementation inconsistencies. The cost of the forward evaluations of the decoder are
    comparable instead, not changing much the online costs. The training of the
    CAEs could be furtherly reduced with transfer
    learning~\cite{goodfellow2016deep} or preprocessing steps that enlight some
    features of the dynamics that are more easily learnable, as was done
    in~\cite{fresca2022pod}. Also the number of training snapshots could be
    optimized furtherly for the test cases presented. The same observations
    apply also for the compressed decoders. Moreover, a parallel
    implementation of the training procedures on more than one GPU is mandatory to achieve
    competitive computational costs.
    \item Simultaneous CAE and compressed decoder/LSTM training. The additional
    costs of the LSTM and compressed decoder training could be cut with a
    unified training of the CAE and compressed decoder: after some epochs, the training of the LSTM or
    compressed decoder could be switched on and performed at the same time of
    the CAE's since the only additional
    information, a part from the restriction of the snapshots into the magic
    points, is the latent dynamics coordinates learned anyway during the CAE's
    optimization.
    \item Increasing the speed-up of the non-linear manifold ROMs. The bottleneck
    for the efficiency of the NM-LSPG-ROC-TS and NM-LSPG-ROC ROMs in the online stage is the
    cost of the compressed decoder or CAE's decoder forward. Regarding the NLC
    model, our results for the
    average time step of the NM-LPSG-GNAT-TS and NM-LSPG-GNAT ROMs from Table~\ref{tab: decoder ts ncl} are comparable if not
    lower than the approximate time step of $7-8$ milliseconds for the NM-LSPG-GNAT with shallow
    autoencoders ROM presented in~\cite{kim2020fast}. The difference is that now
    the FOM implemented with the FVM in OpenFoam~\cite{OF} takes $2.42$ seconds
    for $2000$ time steps
    instead of $140.67$ seconds for $1500$ time steps in~\cite{kim2020fast} with finite differences
    for a smilar test case, i.e. 2d burgers equation instead of our non-linear
    conservation law. To show a more consistent speed-up w.r.t. the FOM, as for
    the SWE test case, the number of degrees of freedom could be increased
    without influencing the NM-LSPG-ROC-TS ROM since it is not dependent on the
    FOM's dimension. In a weaker sense also the NM-LSPG-ROC ROM is also independent
    on the number of degrees of freedom of the FOM, a part from the weights of
    the CAE's decoder that concur in increasing the cost of a single forward in
    the online stage.
    \item Generalization error of LSTM and CAE and additional inductive biases.
    The generalization error of the NN employed depends on a lot of factors, from
    the number of training samples to the regularization term in the loss, the
    architectures, etc. It is thus difficult to predict how much the accuracy of
    the predictions will decay outside the training range. Adding a
    physics-informed term in the loss of the CAE does not provide an improvement
    of the reconstruction error if enough training data are employed as in our
    test cases. Some regularization properties could be imposed in the latent
    dynamics to facilitate the evolution of the NM-LPSG-ROC
    ROMs, for example imposing a linear latent dynamics could be beneficial
    \item Purely data-driven LSTM ROM vs NM-LSPG-ROC and NM-LSPG-ROC-TS ROMs.
    One crucial difference is interpretability: in the first case, the latent
    dynamics is obtained training the LSTM to approximate the latent
    coordinates, in the second case, the latent dynamics is evolved in time
    minimizing the hyper-reduced residual based on the physical model's
    equations. Regarding the LSTM predictions, we have seen in the test cases that despite the higher accuracy
    in the training range and the low computational online cost, there might be
    other issues in the extrapolation regime: in the NLC test case, the
    accuracy was sensitively lower in the extrapolation regime, even if it could
    be improved in theory increasing the layers and nodes of the LSTM in
    exchange for a higher training cost; in the SWE test case the initial time
    step from $0.01$ to $0.017$ seconds could not be well-approximated due to
    different scales and overlappings in the latent dynamics. The NM-LSPG-ROC
    and NM-LSPG-ROC-TS mitigated in some sense these issues, delegating
    less effort in the tuning of the hyper-parameters of the LSTM and increasing
    the interpretability of the results. Moreover the LSTM NN, approximate the
    dynamics only every time step imposed by the training input-output pairs,
    while the non-linear manifold ROMS, can in principle approximate the latent
    dynamics with an arbitrarly small time step, since the decoder provides a
    continuous approximation of the solution manifold and the dynamics is
    evolved based on a numerical scheme with changable time step. With respect
    to purely data-driven methods though, non-linear manifold methods require a
    lower reconstruction error of the CAE since every successive ROMs' dynamics
    evolution depends on how well the intermediate solutions are reconstructed
    through the decoder.
    \item Learn the solution manifolds with autoencoders in unstructured meshes. The natural question of how to extend the CAE architecture to 3D or
    unstructured meshes, is being currently studied. In the literature there are
    already interesting results that employ graph neural networks and their variants to find latent representations
    of 3d simulations~\cite{pfaff2020learning}.
\end{itemize}

\section{Conclusions and perspectives}
\label{sec:conclusions}

We have developed two new hyper-reduced non-linear manifold ROMs: NM-LSPG-ROC and
NM-LSPG-ROC-TS, that can be converted in NM-LSPG-GNAT and NM-LSPG-GNAT-TS. In
NM-LPSG-ROC-TS the
residuals of the NM-LSPG ROM are hyper-reduced with over-collocation, while the
decoder of the CAE is approximated with teacher-student training into a
compressed decoder. In NM-LSPG-ROC only the residuals' hyper-reduction is
carried out. The methods perform similarly in accuracy and computational cost
w.r.t. the NM-LSPG-GNAT with shallow autoencoders ROM introduced
in~\cite{kim2020fast}, for a similar test case. The flexibility of our method permits to change the CAE
architecture depending on the problem at hand, without imposing too many
constraints on its structure, in order to reach the convergence faster and
achieve a lower reconstruction error.

With respect to purely data-driven ROMs built on the CAE's solution manifold,
NM-LSPG-ROC and NM-LSPG-ROC-TS provide more interpretable and, in the extrapolation regimes, sometimes more
accurate predictions, and they need less
hyper-parameters
tuning, once the CAE is trained. Moreover, the methods devolped are
equations-based rather than fully-intrusive, and exploit the physics of the
model to evolve the latent dynamics. It is crucial, though, that the CAE's
reconstruction error is sufficiently low, since the latent dynamics needs to be
computed with numerical schemes that rely on the accuracy of the reconstructed
solutions through the decoder of the CAE or the compressed decoder. We base this
observations on the results obtained for two parametric non-linear time-dependent benchmarks presented
in the numerical results
section~\ref{sec:results}, that is a 2d non-linear conservation law model (NLC)
and a 2d shallow water equations model. Despite the speed-up is not achieved or
not significant with respect to the FOMs, we reached satisfactory results in
terms of the accuracy and the latent or reduced dimension of the ROMs.

Future directions involve the implementation of the developed ROMs in
real applications with higher computational costs and degrees of freedom, such
that a more evident speed-up is reached. This may involve the developement of
autoencoders for more general mesh-based simulations, rather than orthogonal grids.

\appendix
\setcounter{equation}{0}
\renewcommand\theequation{\arabic{equation}}

\section{Neural networks' architectures}
\label{sec:NN}

The architectures of the CAEs are shown in Tables~\ref{tab: cae ncl},~\ref{tab:
cae swe}, of the compressed decoders in Tables~\ref{tab: dec ncl} and \ref{tab:
dec swe}, and of the LSTMs in Table~\ref{tab: lstm}. We mainly employ the Exponential Linear
Unit (ELU) and Rectified Linear Unit (ReLU) activation functions. The padding is
symmetric. The labels \textit{Conv2d},
\textit{ConvTr2d} and \textit{ConvTr2dRec}, stand for 2d convolutions, transposed 2d
convolutions~\cite{NEURIPS2019_9015}, and 2d transposed convolution associated
to a recurrent layer, i.e. they are summed to the previous layer and then
passed to an ELU activation function.

\begin{table}[htp!]
    \centering
    \caption{Nonlinear Conservation Law model's Convolutional Autoencoder.}
    \footnotesize
    \begin{tabular}{ l | c | c | c }
        \hline
        \hline
        Encoder & Activation & Weights & Padding  \\
        \hline
        \hline
        Conv2d & ELU   & [2, 8, 5, 5] & 0 \\
        \hline
        Conv2d & ELU   & [8, 16, 3, 3] & 1 \\
        \hline
        Conv2d & ELU   & [16, 32, 3, 3] & 1 \\
        \hline
        Conv2d & ELU   & [32, 64, 3, 3] & 1 \\
        \hline
        Conv2d & ELU   & [64, 128, 2, 2] & 1 \\
        \hline
        Linear & ELU   & [1152, 4] & - \\
        \hline
        \hline
    \end{tabular}
    \hspace{1mm}
    \begin{tabular}{ l | c | c | c }
        \hline
        \hline
        Decoder & Activation & Weights & Padding  \\
        \hline
        \hline
        Linear & ELU   & [4, 1152] & - \\
        \hline
        ConvTr2d & ELU   & [128, 64, 2, 2] & 1 \\
        \hline
        ConvTr2d & ELU   & [64, 32, 3, 3] & 1 \\
        \hline
        ConvTr2d & ELU   & [32, 16, 4, 4] & 1 \\
        \hline
        ConvTr2d & ELU   & [16, 8, 4, 4] & 0 \\
        \hline
        ConvTr2d & ReLU   & [8, 2, 4, 4] & 1 \\
        \hline
        \hline
    \end{tabular}
    \label{tab: cae ncl}
\end{table}

\begin{table}[htp!]
    \centering
    \caption{Shallow Water Equations model's Convolutional Autoencoder. The decoder for $U$ has two recurrent layers.}
    \footnotesize
    \begin{tabular}{ l | c | c | c }
        \hline
        \hline
        Encoder & Activation & Weights & Padding  \\
        \hline
        \hline
        Conv2d & ELU   & [2, 8, 5, 5] & 0 \\
        \hline
        Conv2d & ELU   & [8, 16, 3, 3] & 1 \\
        \hline
        Conv2d & ELU   & [16, 32, 3, 3] & 1 \\
        \hline
        Conv2d & ELU   & [32, 64, 3, 3] & 1 \\
        \hline
        Conv2d & ELU   & [64, 128, 2, 2] & 1 \\
        \hline
        Linear & ELU   & [1152, 4] & - \\
        \hline
        \hline
    \end{tabular}
    \hspace{1mm}
    \begin{tabular}{ l | c | c | c }
        \hline
        \hline
        Decoder h& Activation & Weights & Padding  \\
        \hline
        \hline
        Linear & ELU   & [4, 1152] & - \\
        \hline
        ConvTr2d & ELU   & [240, 120, 2, 2] & 1 \\
        \hline
        ConvTr2d & ELU   & [120, 60, 3, 3] & 1 \\
        \hline
        ConvTr2d & ELU   & [60, 30, 4, 4] & 1 \\
        \hline
        ConvTr2d & ELU   & [30, 15, 4, 4] & 0 \\
        \hline
        ConvTr2d & ReLU   & [15, 1, 4, 4] & 1 \\
        \hline
        \hline
    \end{tabular}

    \vspace{2 mm}

    \begin{tabular}{ l | c | c | c }
        \hline
        \hline
        Decoder U& Activation & Weights & Padding  \\
        \hline
        \hline
        Linear & ELU   & [4, 2700] & - \\
        \hline
        ConvTr2d & ELU   & [300, 75, 2, 2] & 1 \\
        \hline
        ConvTr2d & ELU   & [150, 75, 3, 3] & 1 \\
        \hline
        ConvTr2dRec & -   & [150, 75, 3, 3] & 1 \\
        \hline
        ConvTr2d & ELU   & [75, 35, 4, 4] & 1 \\
        \hline
        ConvTr2d & ELU   & [35, 20, 4, 4] & 0 \\
        \hline
        ConvTr2dRec & -   & [35, 20, 4, 4] & 0 \\
        \hline
        ConvTr2d & -   & [20, 2, 4, 4] & 1 \\
        \hline
        \hline
    \end{tabular}
    \label{tab: cae swe}
\end{table}

\begin{table}[htp!]
    \centering
    \caption{Nonlinear Conservation Law model's Compressed Decoders.}
    \footnotesize
    \begin{tabular}{ l | c | c | c | c}
        \hline
        \hline
        Magic Points & $1^{st}$ layer weight & $1^{st}$ layer activation & $2^{nd}$ layer weight & $2^{nd}$ layer activation \\
        \hline
        \hline
        50 & [4, 300] & ELU &[300, 278] & ReLU\\
        \hline
        100 & [4, 350] & ELU &[350, 492] & ReLU\\
        \hline
        150 & [4, 400] & ELU &[400, 670] & ReLU\\
        \hline
        \hline
    \end{tabular}
    \label{tab: dec ncl}
\end{table}

\begin{table}[htp!]
    \centering
    \caption{Shallow Water Equations model's Compressed Decoders: height $h$ and velocity $U$ in order.}
    \footnotesize
    \begin{tabular}{ l | c | c | c | c}
        \hline
        \hline
        Magic Points & $1^{st}$ layer weight & $1^{st}$ layer activation & $2^{nd}$ layer weight & $2^{nd}$ layer activation \\
        \hline
        \hline
        25 & [4, 200] & ELU &[200, 238] & ReLU\\
        \hline
        50 & [4, 200] & ELU &[200, 345] & ReLU\\
        \hline
        100 & [4, 300] & ELU &[300, 590] & ReLU\\
        \hline
        150 & [4, 300] & ELU &[300, 654] & ReLU\\
        \hline
        \hline
    \end{tabular}
    \vspace{1 mm}

    \begin{tabular}{ l | c | c | c | c}
        \hline
        \hline
        Magic Points & $1^{st}$ layer weight & $1^{st}$ layer activation & $2^{nd}$ layer weight & $2^{nd}$ layer activation \\
        \hline
        \hline
        25 & [4, 400] & ELU &[400, 478] & -\\
        \hline
        50 & [4, 400] & ELU &[400, 960] & -\\
        \hline
        100 & [4, 600] & ELU &[600, 1180] & -\\
        \hline
        150 & [4, 600] & ELU &[600, 1308] & -\\
        \hline
        \hline
    \end{tabular}
    \label{tab: dec swe}
\end{table}

\begin{table}[htp!]
    \centering
    \caption{Nonlinear Conservation Law model's and Shallow Water Equations model's Long-Shot Term Memory Neural Network.}
    \footnotesize
    \begin{tabular}{ l | c | c | c }
        \hline
        \hline
         & input dim & output dim & number of LSTM layers\\
        \hline
        LSTM layer & 2 & 100 & 2\\
        \hline
        \hline
    \end{tabular}
    \vspace{ 1 mm}

    \begin{tabular}{ l | c | c | c | c}
        \hline
        \hline
         & $1^{st}$ layer weight & $1^{st}$ layer activation & $2^{nd}$ layer weight & $2^{nd}$ layer activation\\
        \hline
        Linear encoding layer & [100, 50] & ELU &[50, 4] & - \\
        \hline
        \hline
    \end{tabular}
    \label{tab: lstm}
\end{table}

\section*{Acknowledgements}
This work was partially funded by European Union Funding for Research and
Innovation --- Horizon 2020 Program --- in the framework of European Research
Council Executive Agency: H2020 ERC CoG 2015 AROMA-CFD project 681447 ``Advanced
Reduced Order Methods with Applications in Computational Fluid Dynamics'' P.I.
Professor Gianluigi Rozza. We also acknowledge the PRIN 2017 “Numerical Analysis for Full and Reduced Order Methods for the efficient and accurate solution of complex
systems governed by Partial Differential Equations” (NA-FROM-PDEs)

\bibliographystyle{abbrv}
\bibliography{biblio}


\end{document}